\documentclass[12pt]{amsart}
\usepackage[margin=1in]{geometry}
\usepackage[utf8]{inputenc}
\usepackage{amsmath}
\usepackage{amsthm}
\usepackage{amssymb}
\usepackage{tikz-cd}
\usepackage{multirow}
\usepackage{lscape}
\usetikzlibrary{decorations.markings}
\tikzstyle{base}=[circle, draw, fill=black,inner sep=0pt, minimum width=4pt]
\tikzstyle{affine}=[circle, draw, fill=red,inner sep=0pt, minimum width=4pt]
\tikzstyle{affine2}=[circle, draw, fill=red,inner sep=0pt, minimum width=8pt]
\tikzstyle{affine3}=[circle, draw, fill=red,inner sep=0pt, minimum width=12pt]
\tikzstyle{affine4}=[circle, draw, fill=red,inner sep=0pt, minimum width=16pt]
\tikzstyle{invis}=[circle,inner sep=0pt, minimum width=4pt]
\tikzstyle{fat2}=[circle,draw,fill=black,inner sep=0pt, minimum width=8pt]
\tikzstyle{fat3}=[circle,draw,fill=black,inner sep=0pt, minimum width=12pt]
\tikzstyle{fat4}=[circle,draw,fill=black,inner sep=0pt, minimum width=16pt]
\usepackage{graphicx}
\usepackage{subfiles}
\usepackage{subcaption}
\usepackage{mathtools}
\usepackage{hyperref}
\usepackage[citation-order]{amsrefs}

\theoremstyle{definition}
\newtheorem{example}{Example}[section]
\newtheorem{definition}[example]{Definition}

\newtheorem{remark}[example]{Remark}
\theoremstyle{plain}
\newtheorem{theorem}[example]{Theorem}

\newtheorem{lemma}[example]{Lemma}

\newcommand{\ClusterAlgebra}[1]{\mathcal{C}(#1)}
\newcommand{\ClusterComplex}[1]{\mathcal{M}_{#1}}
\newcommand{\Dhat}[1]{\widehat{D_{#1}}}
\newcommand{\A}{\mathcal{A}}
\newcommand{\As}{\Tilde{\mathcal{A}}}
\newcommand{\Gr}{Gr}

\title{Special Folding of Quivers and Cluster Algebras}
\author{Dani Kaufman}
\date{April 2023}
\address{University of Copenhagen \\
         Department of Mathematical Sciences \\
         2100 Copenhagen \o, Denmark }
\email{dk@math.ku.dk}
\thanks{The author was supported by the Danish National Research Foundation (CPH-GEOTOP-DNRF151)\\
2020 {\em Mathematics Classification. } Primary 13F60, Secondary 05E40\\
Special thanks to Zachary Greenberg for his programs for cluster algebra computations}

\begin{document}

\begin{abstract}
    We give a precise definition of folded quivers and folded cluster algebras. We give many examples of including some with finite mutation structure that do not have analogues in the unfolded cases. We relate these examples to the finite mutation type quivers $X_6$ and $X_7$. We also construct a folded cluster algebra associated to triangulations of punctured surface which allow for triangulations self-folded triangles. We give a simple construction of a folded cluster algebra for which the cluster complex is a generalized permutohedron. 
\end{abstract}

\maketitle

\section{Introduction}

We consider certain families and examples of \emph{special folding} of quivers, their mutations and associated cluster algebras. Essentially, a folding of a quiver is just a grouping of its nodes into disjoint sets where we then do mutations of each group of nodes together. The associated ``folded cluster algebra'' is generated by assigning the same cluster variable to each node in each group. Folding of quivers is a natural and well understood way to produce cluster algebras with non-skew symmetric exchange matrices from ones with skew symmetric exchange matrices, see \cite{Felikson_unfoldings}.

Our notion of a ``special'' folding simply refers to any folding which cannot be represented by a skew symmetrizable exchange matrix. 

These special foldings were originally motivated from a desire to have a cluster algebra associated to a surface with punctures where there are cluster variables which correspond to arcs enclosing self folded triangles in the triangulation. We construct such a cluster algebra in section \ref{sec:punctured} by taking a two fold covering of our surface ramified over the punctures and folding the associated quiver of the covering surface by the deck transformation. The resulting cluster algebra is essentially the usual cluster algebra associated to a surface, but also contains cluster variables for the self folded arcs. These cluster variables can be written as monomials in the other cluster variables and are actually already contained in the usual surface cluster algebra. 

There are remarkable similarities between our construction and the constructions of generalized cluster algebras associated to surfaces with orbifold singularities of \cite{Felikson_triangulated_orbifolds}. Their folding in the case of a disk with $n+1$ marked points and a $\mathbb{Z}_2$ orbifold point produces a type $B_{n}$ cluster algebra from a type $D_{n+1}$ cluster algebra. However, in our case we cyclically fold a punctured disk with $2n$ marked points to form what we call a type $\Dhat{n}$ cluster algebra from a type $D_{2n}$ algebra. We remark that the cluster combinatorics of our $\Dhat{n}$ cluster algebra and the $B_{n}$ algebra are very similar, but not identical. 

Along the way, we exhibit in section \ref{sec:examples} many other examples of quivers with special foldings that do not come from surfaces. Most of the examples we consider are foldings of mutation infinite quivers where the folded quiver is actually mutation finite. In these cases we use the program written by Z. Greenberg available at \url{https://bitbucket.org/zng42/clusteralgebras/src/master/} to compute the mutation class.
Many, but not all, of these examples are actually the generalized cluster structures related to the cyclic symmetry loci in Grassmannians described by Fraser in \cite{fraser2020cyclic}. 
It would be interesting to wonder if a classification of finite and mutation finite folded quivers would be possible. We find the number of new examples that are possible to construct to be quite surprising, but they do seem to mirror the classification in the unfolded case. 

We finally construct a cluster structure where the cluster complex is the generalized permutohedron in section \ref{sec:permuta}. The cluster variables in this algebra are not very interesting, they can all easily be seen to be monomials in the initial cluster variables.

\section{Cluster Algebras From Quivers}
	
First we recall some notation and constructions related cluster algebras associated to quivers.	
	
\subsection{Quivers and mutations}
    
    \begin{definition}
    A \emph{quiver} is a finite directed graph without self loops or 2-cycles. 
    \end{definition}
    
    Such a quiver is equivalent to a skew symmetric matrix, usually called the exchange matrix in the context of cluster algebras. We will prefer to work with quivers since it allows for easier graphical representation. Let $Q$ be a quiver.
    \begin{definition}
    \emph{Mutation} of $Q$ at node $k$ is an operation which produces a new quiver, $\mu_k(Q)$, where the arrows of $Q$ change by the following rule:
    \begin{enumerate}
        \item For every pair of arrows $(i \rightarrow k), (k \rightarrow j)$ add a arrow from node $i$ to node $j$, cancelling any 2-cycles between nodes $i$ and $j$. 
        \item Reverse every arrow incident to $k$. 
    \end{enumerate}
    \end{definition}
    
    \begin{definition}
    We call two quivers \emph{mutation equivalent} if there is a sequence of mutations that transforms one into the other. We call the set of quivers mutation equivalent to a given quiver its \emph{mutation class}. Quivers with a finite mutation class are called \emph{mutation finite}.
    \end{definition}
    
    When we name a quiver, we usually are giving a name to its mutation class.
    
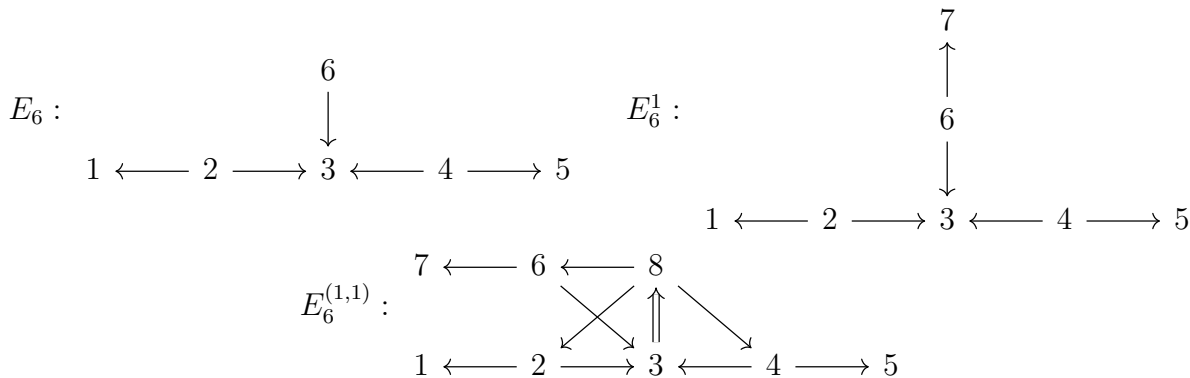
\begin{figure}
        \centering
        $E_6:$
        \begin{tikzcd}
  &                       & 6 \arrow[d] &                       &   \\
1 & 2 \arrow[l] \arrow[r] & 3           & 4 \arrow[l] \arrow[r] & 5
\end{tikzcd} \quad
$E_6^1:$ 
\begin{tikzcd}
  &                       & 7                     &                       &   \\
  &                       & 6 \arrow[d] \arrow[u] &                       &   \\
1 & 2 \arrow[l] \arrow[r] & 3                     & 4 \arrow[l] \arrow[r] & 5
\end{tikzcd}
$E_6^{(1,1)}:$
\begin{tikzcd}
7 & 6 \arrow[rd] \arrow[l] & 8 \arrow[l] \arrow[ld] \arrow[rd] &                       &   \\
1 & 2 \arrow[l] \arrow[r]  & 3 \arrow[u, Rightarrow]           & 4 \arrow[l] \arrow[r] & 5
\end{tikzcd}
        \caption{Finite, Affine, and Elliptic quivers of types $E_6$}
        \label{fig:e6}
    \end{figure}

    The classification of mutation finite quivers is well understood, see \cite{felikson_finitemutations}. These quivers either come from triangulations of surfaces or are one of 11 exceptions. 9 of the exceptions are related to type $E_n$ Dynkin diagrams; each $E_n$ gives a finite, affine, and elliptic Dynkin diagram, see figure \ref{fig:e6}. The other two exceptions are simply called $X_6$ and $X_7$ with 6 and 7 nodes respectively.
    
 \subsection{Folding quivers}   
    
    The relationship between the classical folding of simply laced Dynkin diagrams to form non simply laced diagrams can be extended to quivers. 
    
    Essentially, to fold a quiver we will group its nodes into disjoint sets and do mutations by requiring that we mutate all the nodes of a given set together. We call the operation of mutating each of the elements of a set of nodes in turn a ``group mutation''. 
    
    In general a group mutation will depend on the order the nodes are mutated. However if all the mutations commute, the order doesn't matter. A sufficient condition for all the mutations to commute is that there are no arrows between any two nodes in the same group. 
    
    With this in mind, we have the following definition:
    \begin{definition}
    A \emph{folding} of a quiver, $Q$, with $n$ nodes is a choice of $k$ non empty and disjoint sets of nodes whose union contains all of the nodes of $Q$ satisfying the following conditions called the \emph{folding conditions}:
    \begin{enumerate}
        \item The nodes contained in a given set have no arrows between themselves.
        \item Condition 1 is satisfied after any number of group mutations of these fixed sets.
    \end{enumerate}
    We call a folding \emph{valid} if it satisfies the folding conditions. If a folding is valid, then it defines a \emph{folded mutation class} which is simply the set of folded quivers group mutation equivalent to it.
    \end{definition}
    
    When we talk of a folding of a quiver we are essentially creating a new quiver where we think of each group as a single node. When we name a folded quiver we are generally referring to its folded mutation class.
    
    \begin{remark}
    There is a natural generalisation of the folding conditions to allow for foldings where there are arrows between the nodes in the same group. If the subquiver of nodes in a group has a "red to green" mutation sequence that is order two (up to an automorphism of the subquiver), then applying this transformation should behave like a generalized mutation of the group. Our folding condition would then require that this is always possible. 
    \end{remark}
    
    These properties will always be satisfied when $Q$ is the quiver of a simply laced finite or affine Dynkin diagram and the sets are given by the orbits of nodes under an automorphism of $Q$. In this situation the group mutations will describe the usual matrix mutations associated to a ``skew symmetrizable'' matrix, see \cites{Calg:1,Felikson_triangulated_orbifolds} 
    
    \begin{definition}
    We call a folding which can be represented by a skew symmetrizable matrix a \emph{Standard Folding}. Any other folding is referred to as a \emph{Special Folding}.
    \end{definition}
    
    It is often the case that we want to fold a quiver by the action of an automorphism or subgroup of automorphisms. 
    
    \begin{lemma}\label{lem:automorpic_mutations}
    Let $\tau$ be an automorphism of a quiver $Q$ and let $O$ be an orbit of nodes under the action of $\tau$ with no arrows amongst the nodes in $O$. Then a group mutation at $O$ produces a new quiver on which $\tau$ acts via an automorphism.   
    \end{lemma}
    \begin{proof}
    We need to check that each arrow that is changed or added by the mutation operation is done in a $\tau$ symmetric way. Let $(i \rightarrow k), (k \rightarrow j)$ be a pair of arrows with $k \in O$ and $i,j \notin O$. Then since we mutate at all nodes in $O$ we add the same number of arrows $\tau^n(i) \to \tau^n(j)$ for all $n$, which preserves $\tau$ symmetry. The second step of mutation also clearly preserves $\tau$ symmetry.
    \end{proof}
    
    It may be the case, however, that this folding may not actually be a valid folding. The next lemma will give us a large class of examples  
    
    \begin{lemma}\label{thm:order_2folds}
    If a quiver $Q$ has an order 2 automorphism, $\tau$, then the folding of $Q$ by the orbits of $\tau$ is valid. 
    \end{lemma}
    \begin{proof}
    First we note that the orbits of $\tau$ each contain one or two nodes. We just need to make sure that the orbits with two nodes never end up with an arrow between them after some number of group mutations. To see this, we notice that any group with two nodes cannot have an arrow between them since $i,j$ are an orbit then an arrow $(i \to j)$ implies that there must be an arrow $(\tau(i) \to \tau(j))$ which would generate a 2 cycle in $Q$. 
    Now by the previous lemma, we see that group mutations of this folding generate new quivers on which $\tau$ is a automorphism. 
    \end{proof}

\subsection{Cluster algebras from quivers}
    
    Cluster algebras (of geometric type) are formed by starting with an initial seed consisting of a quiver with variables associated to its nodes and applying all possible mutations. 
    
    Let $Q$ be a quiver with $n$ nodes and let $\boldsymbol{z} =\{z_1, \dots,z_n\}$ be algebraically independent elements of $\mathbb{Z}(x_1,\dots,x_n)$.
    \begin{definition}
     A \emph{seed}, $\boldsymbol{i}=(Q,\boldsymbol{z})$ is a pair of a quiver $Q$ and the set of variables $\boldsymbol{z} = \{z_1, \dots,z_n\}$ with each variable associated to a distinct node of $Q$. The set $\boldsymbol{z}$ is called a \emph{cluster} and the variables are called \emph{cluster variables}. 
    \end{definition}
    
    We can mutate a seed at a node $k$ by mutating $Q$ as before and replacing the variable $z_{k}$ associated with node $k$ with a new variable $z_{k}'$ satisfying
    \begin{equation}
        z_{k}\cdot z_{k}' = \prod_{(i \rightarrow k)} z_{i} + \prod_{(k \rightarrow j)} z_{j}.
    \end{equation} 
    This procedure generates a new seed, and a new cluster of variables. 
    
    We may also mark some nodes of $Q$ as frozen, and do not allow mutations at these nodes. The variables associated to these nodes are called "frozen variables" and will appear in every cluster. 
    
    \begin{definition}
    The $\mathbb{Z}-$algebra generated by all of the cluster variables obtained from all possible mutations of a seed is the \emph{cluster algebra} associated with that seed, denoted by $\ClusterAlgebra{Q}$ where $Q$ is the quiver underlying the initial seed.
    \end{definition}
    
    The primary examples of cluster algebras we will reference are the Grassmannian cluster algebras denoted by $\Gr(n,k)$ defined by Scott \cite{scott_grassmannains}.

\subsection{The cluster complex}

    While we wont need many of the details, we will quickly recall the notion of the cluster complex associated to a cluster algebra.
    
    For any cluster algebra there is an associated simplicial complex $\ClusterComplex{Q}$ called the cluster complex. This complex is defined in detail in \cite{Calg:1,Fock_gonch}.
    
    \begin{definition}
    Two cluster variables are \emph{compatible} if they appear in a cluster together. The cluster complex is the clique complex of the compatibility rule for cluster variables.
    \end{definition}
    
    In other words, the $k$-dimensional simplices of $\ClusterComplex{Q}$ correspond to size $k$ collections of mutually compatible cluster variables in $\ClusterAlgebra{Q}$.  Vertices correspond to an individual cluster variables and each edges connects two cluster variables when they can be found in a cluster together. The maximal dimension simplices correspond to the clusters of $\ClusterAlgebra{Q}$.
    
    \begin{definition}
    The dual complex of the cluster complex is called the \emph{exchange complex}
        The 1-skeleton of the exchange complex is called the \emph{exchange graph} of the cluster algebra. 
    \end{definition}
    
    The vertices of this graph correspond to clusters and the edges correspond to mutations between clusters. 
    
    The group of symmetries of these complexes is often of interest and is called the ``cluster modular group''.

\subsection{Folding Cluster algebras}

If a quiver satisfies the folding condition, we are interested in whether there is an associated folded cluster algebra akin to the usual cluster algebras associated with non-skew symmetric exchange matrices.
\begin{definition}
A folded quiver satisfies the \emph{folded cluster condition} or simply called a \emph{Cluster folding} if the initial seed where each group is assigned the same cluster variable for all its nodes retains this property after any number of group mutations.  We call the cluster algebra generated this way the \emph{folded cluster algebra} associated to this folded quiver.
\end{definition}

\begin{lemma}\label{thm:folded_clusters}
If the folding of $Q$ under an automorphism $\tau$ is valid, then it satisfies the folded cluster condition.
\end{lemma}

\begin{proof}
We just need to check that each group mutation produces the same cluster variable on each node in a given group. This follows simply by seeing that $\tau$ acts on each exchange relation, and that the action of $\tau$ on a given initial cluster variable is trivial since we start with the same variable on each orbit.
\end{proof}

\begin{remark}
We could generalize this definition by instead of placing the same variable on each node in a group, but instead that some relationship is always satisfied between the variables on each node in a group.
\end{remark}

Even if a folding is valid but not cluster, we can still look at the subalgebra of the cluster algebra associated with the unfolded quiver which is picked out by the group mutations. In this way we can still look at an exchange complex for a non-cluster folding.

\subsection{Duality of foldings}

There is a natural generalization of notion of Langlands Duality displayed by the the type $B_n$ and $C_n$ cluster algebras.
    
Let $Q,R$ be two quivers with valid foldings in groups $G_1, \dots ,G_k$ and $H_1, \dots H_k$.
    \begin{definition}
We call these foldings \emph{Langlands Dual} if there is a cellular isomorphism, $\phi$, between their cluster complexes which gives correspondence between the groups of each folding such that the ratios $|G_i|/|G_j| = |\phi(G_j)|/|\phi(G_i)|$.
\end{definition}

We will use this notion to relate some new foldings with other previously known quivers.

\section{Examples of Special Folding With Finite Cluster Complexes}\label{sec:examples}

In this section we represent foldings with $n$ groups by folding nodes $1,2,\dots,n$ with other nodes of the quiver so that each group is identified by a node between 1 and $n$. In the figures the groups will usually be shown vertically. We show pictures of the exchange complexes of many of the algebras mentioned here in the appendix \ref{sec:appendix}.

\subsection{Finite cluster foldings}

Our first example of a special folding is the basic example from which most of our study will continue. Let $D_n^\circ$ be the quiver with $n$ nodes connected in a oriented cycle, with the nodes numbered in increasing order.

The quiver $D_{2n}^\circ$ has a special folding with $n$ groups where node $i$ is grouped with node $i+n \mod 2n$. We will call this folded quiver $\Dhat{n}$. This is a valid folding which satisfies the folded cluster condition by lemmas \ref{thm:folded_clusters} and \ref{thm:order_2folds}

\begin{figure}
    \centering
\begin{tikzcd}
1 \arrow[r]   & 2 \arrow[ld] \arrow[r]   & 3   & \dots \arrow[l] \arrow[r] & n  \\
n+1 \arrow[r] & n+2 \arrow[lu] \arrow[r] & n+3 & \dots \arrow[l] \arrow[r] & 2n
\end{tikzcd}\\
\vspace{.5cm}
\begin{tikzcd}
1 \arrow[r, bend right] & 2 \arrow[r] \arrow[l, bend right] & 3 & \dots \arrow[l] \arrow[r] & n
\end{tikzcd}
    \caption{The representative quiver for $\Dhat{n}$.}
    \label{fig:dn_hat_dynk}
\end{figure}
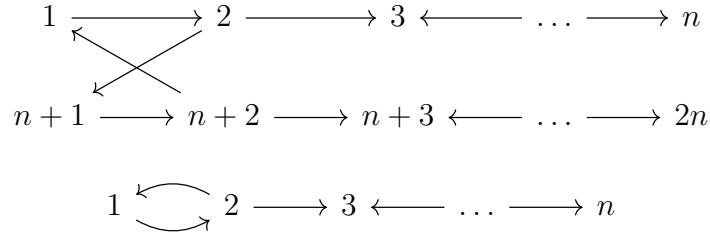

The $\Dhat{n}$ quiver is mutation equivalent to the first quiver shown in figure \ref{fig:dn_hat_dynk}. It is interesting to note that the  exchange relations generated from the group mutations starting from the folding of $D_{2n}^\circ$ are exactly the same as those generated from a quiver which is an oriented $n-$cycle. Once we mutate to a quiver like the one shown in figure \ref{fig:dn_hat_dynk}, the exchange relations are different. This shows that $\Dhat{n}$ cannot be represented by a skew-symmetrizable matrix.  

Many examples of special foldings will involve the subquiver of the first two groups of this quiver, so as a short hand we represent this particular folding (when the folding in unambiguous) by the quiver with a two cycle also shown in figure \ref{fig:dn_hat_dynk}. 


\begin{figure}
    \centering
    $\widehat{E_6}:$ \begin{tikzcd}
1 \arrow[r] & 2 \arrow[r, bend right] & 3 \arrow[r] \arrow[l, bend right] & 4
\end{tikzcd} \quad
$\widehat{E_7}:$
\begin{tikzcd}
1 \arrow[r] & 2 \arrow[r, bend right] & 3 \arrow[r] \arrow[l, bend right] & 4 & 5 \arrow[l]
\end{tikzcd}\\ 
\vspace{.5cm}
$\widehat{E_8}:$ 
\begin{tikzcd}
1 \arrow[r] & 2 \arrow[r, bend right] & 3 \arrow[r] \arrow[l, bend right] & 4 & 5 \arrow[l] \arrow[r] & 6
\end{tikzcd}
    \caption{$\widehat{E_n}$}
    \label{fig:E_specials}
\end{figure}
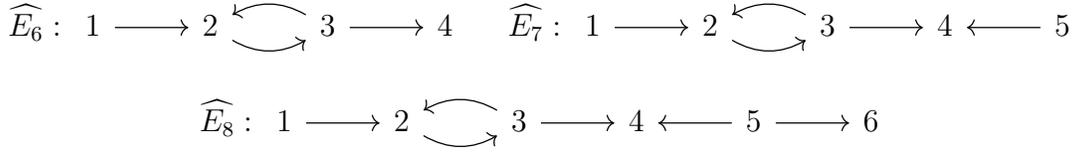

We can also find more quivers with special foldings in the same general family, shown in figure \ref{fig:E_specials}. The folded cluster algebras associated to these quivers should have combinatorics related to the $E_n$ root systems as labeled, but we will not explore these details here.

Interestingly, the folded quivers are finite type for each of these examples, but the unfolded quivers are finite, finite mutation type, and infinite types respectively. This implies that the mutation type of a quiver does not need to be preserved after folding. 
The unfolded quivers of the $\widehat{E_n}$ quivers are mutation equivalent to those of the Grassmannian cluster algebras $\Gr(3,8),\Gr(3,9),\Gr(3,10)$ respectively. 

The correspondence between the usual $D$ and $E$ Dynkin diagrams and these special quivers can be seen by replacing the 2 nodes with a 2 cycle with 4 nodes connected as a $D_4$ diagram. Interestingly, the 2 cycle unfolds to a 4 cycle which is mutation equivalent to the $D_4$ diagram.

\subsection{Examples of non cluster foldings}

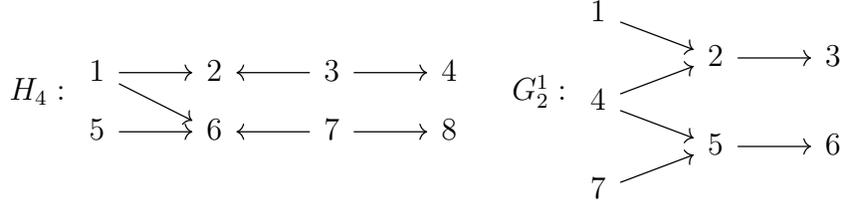
\begin{figure}
    \centering
$H_4:$
        \begin{tikzcd}[row sep =.5em]
1 \arrow[r] \arrow[rd] & 2 & 3 \arrow[r] \arrow[l] & 4 \\
5 \arrow[r]            & 6 & 7 \arrow[l] \arrow[r] & 8
\end{tikzcd}
\quad
$G_2^1:$
    \begin{tikzcd}[row sep=.02em]
1 \arrow[rd]            &             &   \\
                        & 2 \arrow[r] & 3 \\
4 \arrow[ru] \arrow[rd] &             &   \\
                        & 5 \arrow[r] & 6 \\
7 \arrow[ru]            &             &  
\end{tikzcd}
    \caption{Non-cluster foldings from Dynkin Diagrams.}
    \label{fig:h4}
\end{figure}

Our first examples of non-cluster foldings come from folding Dynkin Diagrams. Two examples are shown in figure \ref{fig:h4}. The first comes from the usual relation between the $E_8$ and $H_4$ Coxeter groups. The exchange complex associated to this folding generates the $H_4$ generalized associahedron as constructed in \cite{Reading_generalized_ass}.
The second is a folding of an affine $D_6$ quiver resulting in a quiver with mutation type affine $G_2$. However, there is not a special folded cluster algebra associated with this folding, even though there is an affine $G_2$ algebra. 
 
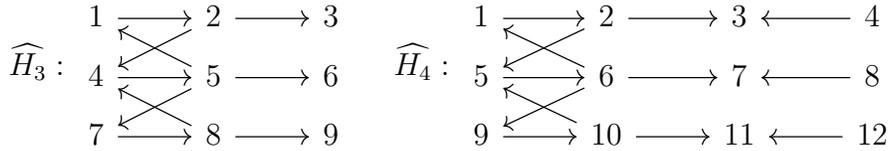
\begin{figure}
    \centering
    $\widehat{H_3}:$
\begin{tikzcd}[row sep=.5em]
1 \arrow[r] & 2 \arrow[ld] \arrow[r]            & 3 \\
4 \arrow[r] & 5 \arrow[lu] \arrow[ld] \arrow[r] & 6 \\
7 \arrow[r] & 8 \arrow[lu] \arrow[r]            & 9
\end{tikzcd}
\quad
$\widehat{H_4}:$
\begin{tikzcd}[row sep=.5em]
1 \arrow[r] & 2 \arrow[ld] \arrow[r]            & 3 & 4 \arrow[l] \\
5 \arrow[r] & 6 \arrow[lu] \arrow[ld] \arrow[r] & 7 & 8 \arrow[l] \\
9 \arrow[r] & 10 \arrow[lu] \arrow[r]            & 11 & 12 \arrow[l]
\end{tikzcd}
    \caption{Special non-cluster folding examples.}
    \label{fig:Strange}
\end{figure}

We can also find some examples which are entirely new. 
The quivers shown in \ref{fig:Strange} have special foldings with the groups shown vertically. These examples also do not satisfy the folded cluster condition. 

The exchange complex of $\widehat{H_3}$ has 32 vertices, 48 edges, 18 faces (4 squares, 8 pentagons, 2 hexagons, and 4 heptagons). We note that this complex has the same face vector as the usual $H_3$ generalized associahedron, see \cite{Reading_generalized_ass}, but it is not exactly the same complex. The exchange complex of $\widehat{H_4}$ is however very different from that of the $H_4$ generalized associahedron.

We tabulate some of the properties of the finite special folded algebras in table \ref{tab:special_ADE}. This includes the $\widehat{H_3}$ and $\widehat{H_4}$ quivers, even though they do not have folded algebras. The total number of clusters in the $\widehat{D_n}$ algebras is computed in section \ref{sec:punctured}. 

\begin{remark}
Like the $H_3$ and $\widehat{H_3}$ examples, the $\Dhat{n}$ exchange complex seems to have the same face vector as the $B_n$ generalized associahedron. We only prove this for the maximal simplicies and 0 simplicies (number of clusters and cluster variables) though. 
\end{remark}

\begin{table}[]
    \centering
    \begin{tabular}{c|c|c|c|c}
       Type  & Unfolded type & Number of cluster variables & Number of clusters & symmetry   \\ \hline
         $\Dhat{n}$ & $D_{2n}$ & $ (n+1)(n)$ & $\binom{2n}{n}$ & $n$   \\
        $\widehat{E_6}$ & $E_8 = \Gr(3,8)$ & 24 & 88 & $4$ \\
        $\widehat{E_7}$ & $E_8^{(1,1)}=\Gr(3,9)  $ & 45 & 432 & $6$ \\
        $\widehat{E_8}$ & $\Gr(3,10)$ & 90 & 2600 & $10$ \\ \hline
        $\widehat{H_3}$ & $E_7^{(1,1)}=\Gr(4,8)$ & 18 & 32 & $4$ \\
        $\widehat{H_4}$ & $\Gr(4,9)$ & 62 & 234 & $9$ \\
     \end{tabular}
    \caption{Cluster combinatorics of the finite type special foldings.}
    \label{tab:special_ADE}
\end{table}

\subsection{Finite mutation type foldings}

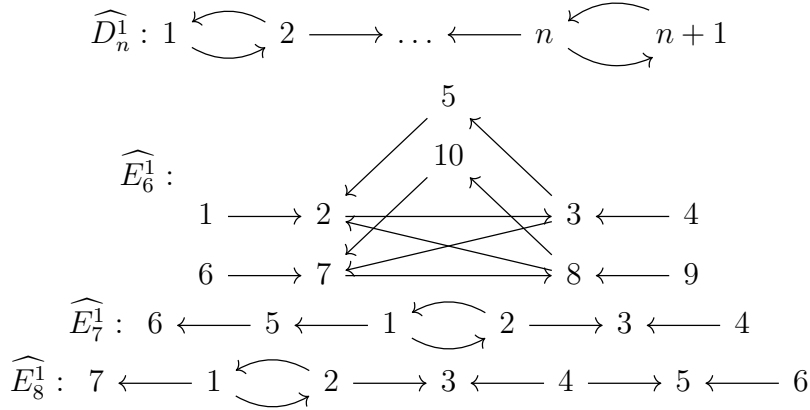
\begin{figure}
    \centering
    $\widehat{D^1_n}:$\begin{tikzcd}
1 \arrow[r, bend right] & 2 \arrow[r] \arrow[l, bend right] & \dots & n \arrow[l] \arrow[r, bend right] & n+1 \arrow[l, bend right]
\end{tikzcd} \\
    $\widehat{E^1_6}:$
\begin{tikzcd}[row sep=.5em]
            &              & 5 \arrow[ldd]  &                           &             \\
            &              & 10 \arrow[ldd] &                           &             \\
1 \arrow[r] & 2 \arrow[rr] &                & 3 \arrow[luu] \arrow[lld] & 4 \arrow[l] \\
6 \arrow[r] & 7 \arrow[rr] &                & 8 \arrow[luu] \arrow[llu] & 9 \arrow[l]
\end{tikzcd} \\
$\widehat{E^1_7}:$
\begin{tikzcd}
6 & 5 \arrow[l] & 1 \arrow[r, bend right] \arrow[l] & 2 \arrow[r] \arrow[l, bend right] & 3 & 4 \arrow[l]
\end{tikzcd} \\
$\widehat{E^1_8}:$
\begin{tikzcd}
7 & 1 \arrow[r, bend right] \arrow[l] & 2 \arrow[r] \arrow[l, bend right] & 3 & 4 \arrow[l] \arrow[r] & 5 & 6 \arrow[l]
\end{tikzcd}
    \caption{Special affine types.}
    \label{fig:E_affines}
\end{figure}

Most of the examples thus far considered produce finitely many clusters. We can construct many new examples of quivers with special foldings where the folding is finite mutation type. For each of the new exceptional example we tabulate the number of quivers in its mutation class in table \ref{tab:my_label}. The first examples are special foldings corresponding to the affine $D$ and $E$ types as shown in figure \ref{fig:E_affines}.

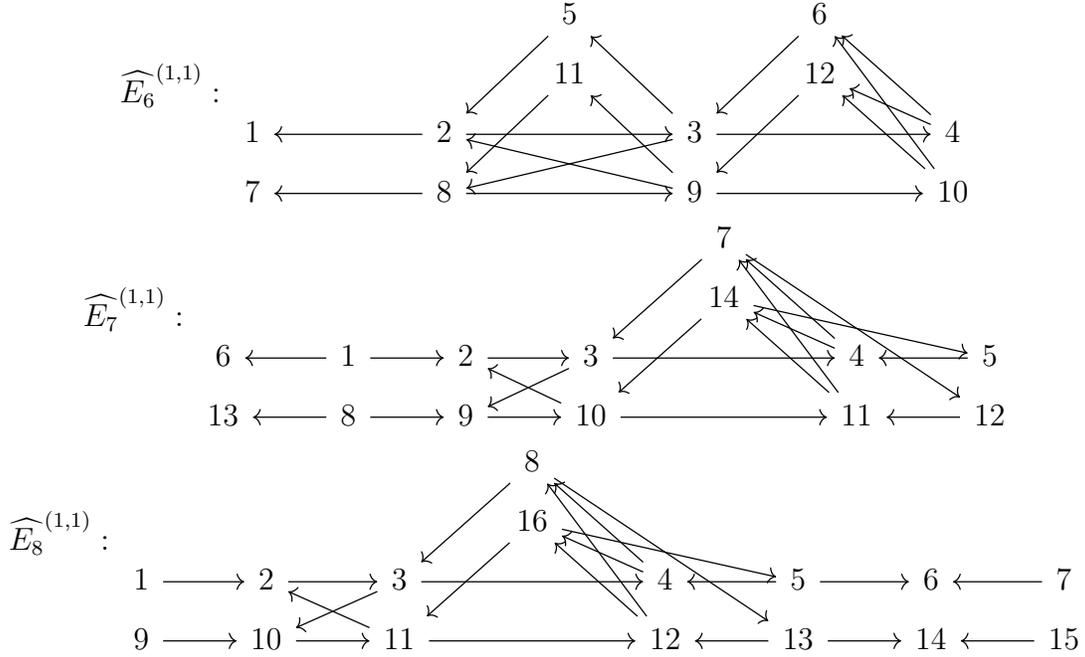
\begin{figure}
    \centering
    $\widehat{E_6}^{(1,1)}:$
\begin{tikzcd}[row sep=.5em]
  &  &                         & 5 \arrow[ldd]  &                                      & 6 \arrow[ldd]  &                             \\
  &  &                         & 11 \arrow[ldd] &                                      & 12 \arrow[ldd] &                             \\
1 &  & 2 \arrow[rr] \arrow[ll] &                & 3 \arrow[lld] \arrow[rr] \arrow[luu] &                & 4 \arrow[luu] \arrow[lu]    \\
7 &  & 8 \arrow[rr] \arrow[ll] &                & 9 \arrow[llu] \arrow[rr] \arrow[luu] &                & 10 \arrow[luu] \arrow[luuu]
\end{tikzcd}
\\

    $\widehat{E_7}^{(1,1)}:$
\begin{tikzcd}[row sep=.5em]
   &                       &             &                          & 7 \arrow[ldd] \arrow[rrddd] &                             &              \\
   &                       &             &                          & 14 \arrow[rrd] \arrow[ldd]  &                             &              \\
6  & 1 \arrow[r] \arrow[l] & 2 \arrow[r] & 3 \arrow[ld] \arrow[rr]  &                             & 4 \arrow[luu] \arrow[lu]    & 5 \arrow[l]  \\
13 & 8 \arrow[r] \arrow[l] & 9 \arrow[r] & 10 \arrow[lu] \arrow[rr] &                             & 11 \arrow[luu] \arrow[luuu] & 12 \arrow[l]
\end{tikzcd}
\\

$\widehat{E_8}^{(1,1)}:$
\begin{tikzcd}[row sep=.5em]
            &              &                          & 8 \arrow[ldd] \arrow[rrddd] &                             &                        &    &              \\
            &              &                          & 16 \arrow[rrd] \arrow[ldd]  &                             &                        &    &              \\
1 \arrow[r] & 2 \arrow[r]  & 3 \arrow[ld] \arrow[rr]  &                             & 4 \arrow[luu] \arrow[lu]    & 5 \arrow[l] \arrow[r]  & 6  & 7 \arrow[l]  \\
9 \arrow[r] & 10 \arrow[r] & 11 \arrow[lu] \arrow[rr] &                             & 12 \arrow[luu] \arrow[luuu] & 13 \arrow[l] \arrow[r] & 14 & 15 \arrow[l]
\end{tikzcd}
    \caption{Special elliptic types.}
    \label{fig:double_extended}
\end{figure}

There are several more quivers with special foldings extending the special Affine E types in a similar way that the usual elliptic Dynkin diagrams of types $E_6^{(1,1)},E_7^{(1,1)},E_8^{(1,1)}$ extend the affine types. These are shown in figure \ref{fig:double_extended}.

\begin{figure}
    \centering
   $\widehat{E_6}^{(2,2)}:$ 
    \begin{tikzcd}[row sep=.5em]
                           & 8 \arrow[rdd]  &                                      & 5 \arrow[ldd]  &                                       & 7 \arrow[ldd]  &                             \\
                           & 16 \arrow[rdd] &                                      & 13 \arrow[ldd] &                                       & 15 \arrow[ldd] &                             \\
1 \arrow[ru] \arrow[ruu]   &                & 2 \arrow[rr] \arrow[ll] \arrow[rddd] &                & 3 \arrow[lld] \arrow[rr] \arrow[luu]  &                & 4 \arrow[luu] \arrow[lu]    \\
9 \arrow[ruuu] \arrow[ruu] &                & 10 \arrow[rr] \arrow[ll] \arrow[rd]  &                & 11 \arrow[llu] \arrow[rr] \arrow[luu] &                & 12 \arrow[luu] \arrow[luuu] \\
                           &                &                                      & 14 \arrow[ruu] &                                       &                &                             \\
                           &                &                                      & 6 \arrow[ruu]  &                                       &                &                            
\end{tikzcd}
\\
    
    $\widehat{E_7}^{(1,2)}:$
\begin{tikzcd}[row sep=.5em]
             &                            & 8 \arrow[rdd] \arrow[llddd] &                         &                          & 7 \arrow[ldd] \arrow[rrddd] &                             &              \\
             &                            & 16 \arrow[rdd] \arrow[lld]  &                         &                          & 15 \arrow[rrd] \arrow[ldd]  &                             &              \\
6 \arrow[r]  & 1 \arrow[ru] \arrow[ruu]   &                             & 2 \arrow[r] \arrow[ll]  & 3 \arrow[ld] \arrow[rr]  &                             & 4 \arrow[luu] \arrow[lu]    & 5 \arrow[l]  \\
14 \arrow[r] & 9 \arrow[ruuu] \arrow[ruu] &                             & 10 \arrow[r] \arrow[ll] & 11 \arrow[lu] \arrow[rr] &                             & 12 \arrow[luu] \arrow[luuu] & 13 \arrow[l]
\end{tikzcd}
    \caption{The further extended quivers. }
    \label{fig:Futher_extended}
\end{figure}
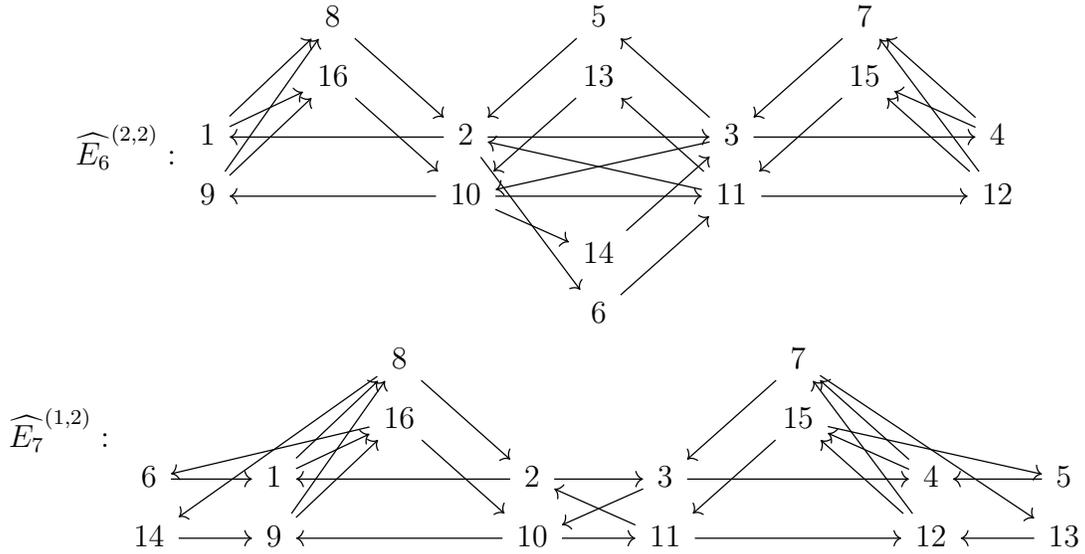

Surprisingly, it is possible to two more groups to the $\widehat{E_6}^{(1,1)}$ quiver and one more group to the $\widehat{E_7}^{(1,1)}$ quiver and maintain finite mutation types. These quivers are shown in figure \ref{fig:Futher_extended}. We also can define $\widehat{E_6}^{(1,2)}$ and $\widehat{E_6}^{(0,2)}$ by freezing groups $8$ and $8,7$ of $\widehat{E_6}^{(2,2)}$ respectively. 

\begin{remark}
The unfolded quivers of $\widehat{E_7}^{(1,2)}$ and $\widehat{E_8}^{(1,1)}$ are mutation equivalent to that of the Grassmannians $\Gr(5,10)$ and $\Gr(3,12)$ respectively. This is a generalization of the fact that the usual  $E_7^{(1,1)}$ and $E_8^{(1,1)}$ algebras are the types of $\Gr(4,8)$ and $\Gr(3,9)$ respectively.
We suspect that the unfolded type of $\widehat{E_6}^{(2,2)}$ has an interesting algebraic interpretation, in analogy to the usual $E_6^{(1,1)}$ being related to the weighted projective space $\mathbb{P}[3,3,3]$.
\end{remark}
 
\subsection{Relation to cyclic symmetry loci of Grassmannians}
jMany of the examples of quivers whose unfolding is a Grassmannian, $\Gr(k,n)$ are associated to a two fold cyclic symmetry. $\Gr(k,n)$ can be thought of as a collection of $n$ $k-$vectors. The $n-fold$ cyclic symmetry of the $n$ vectors extends to a symmetry of the cluster structure. Fraser \cite{fraser2020cyclic} describes the structure of the fixed point set in these cases.  When $n$ is even, we can try to fold the cluster structure by the order two subgroup of this cyclic symmetry. The foldings $\widehat{E_6}, \widehat{E_8}, \widehat{E_7}^{(1,2)} , \widehat{E_8}^{(1,1)} $ all arise this way. The other examples which unfold to Grassmannians do not.

\subsection{Relation to the $X_7$ quiver }

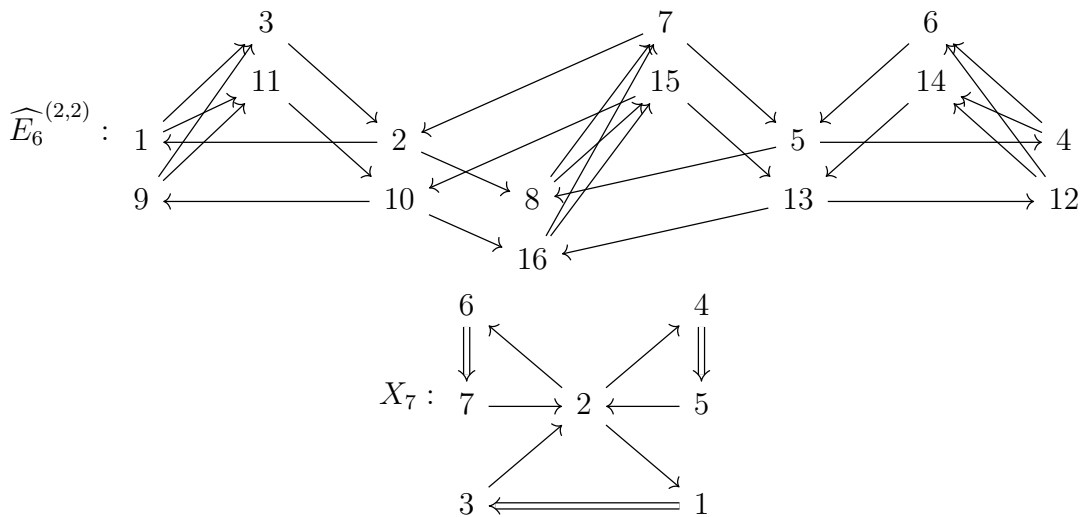
\begin{figure}
    \centering
    $\widehat{E_6}^{(2,2)}:$
\begin{tikzcd}[row sep=.5em]
                           & 3 \arrow[rdd]  &                          &                               & 7 \arrow[lldd] \arrow[rdd]  &                           & 6 \arrow[ldd]  &                             \\
                           & 11 \arrow[rdd] &                          &                               & 15 \arrow[lldd] \arrow[rdd] &                           & 14 \arrow[ldd] &                             \\
1 \arrow[ru] \arrow[ruu]   &                & 2 \arrow[ll] \arrow[rd]  &                               &                             & 5 \arrow[rr] \arrow[lld]  &                & 4 \arrow[luu] \arrow[lu]    \\
9 \arrow[ruu] \arrow[ruuu] &                & 10 \arrow[ll] \arrow[rd] & 8 \arrow[ruu] \arrow[ruuu]    &                             & 13 \arrow[rr] \arrow[lld] &                & 12 \arrow[luu] \arrow[luuu] \\
                           &                &                          & 16 \arrow[ruuu] \arrow[ruuuu] &                             &                           &                &                            
\end{tikzcd}

$X_7:$
\begin{tikzcd}
6 \arrow[d, Rightarrow] &                                    & 4 \arrow[d, Rightarrow]  \\
7 \arrow[r]             & 2 \arrow[lu] \arrow[ru] \arrow[rd] & 5 \arrow[l]              \\
3 \arrow[ru]            &                                    & 1 \arrow[ll, Rightarrow]
\end{tikzcd}
    \caption{A quiver group mutation equivalent to $\widehat{E_6}^{(2,2)}$ and the quiver $X_7$ showing the duality.} 
    \label{fig:Duality_x7}
\end{figure}

The quivers $\widehat{E_6}^{(1,1)}$ and $\widehat{E_6}^{(2,2)}$ are related to the special mutation finite quivers $X_6$ and $X_7$ in the following way: The folding of $X_7$ by grouping nodes $4,6$ and $5,7$ is Langlands dual to a folding of folding of quivers in the group mutation class of $\widehat{E_6}^{(2,2)}$, see figure \ref{fig:Duality_x7}. The first quiver is group mutation to $\widehat{E_6}^{(2,2)}$. We can fold it further by grouping the $1,2,3$ groups with the $4,5,6$ groups respectively. This extra folding is dual to the folding of $X_7$. By freezing node 3 in $X_7$ one sees that there is a folding of $\widehat{E_6}^{(1,1)}$ which is dual to a folding of $X_6$.

\begin{remark}
This shows a natural example where special folding arise. It is not possible to make a quiver which is Langlands dual to a folding of $X_7$ without doubling all of the groups.
\end{remark}

\begin{figure}
    \centering
$Y_5:$
\begin{tikzcd}[column sep=.5em,row sep=.75em]
                &                 &  & 1 \arrow[dddd] \arrow[rdddd]                & 6 \arrow[dddd] \arrow[ldddd]                &  &                  &                   \\
                &                 &  &                                             &                                             &  & 4 \arrow[lllu]   & 9 \arrow[lllu]    \\
3 \arrow[rrruu] & 8 \arrow[rrruu] &  &                                             &                                             &  &                  &                   \\
                &                 &  &                                             &                                             &  & 5 \arrow[llluuu] & 10 \arrow[llluuu] \\
                &                 &  & 2 \arrow[llluu] \arrow[rrru] \arrow[rrruuu] & 7 \arrow[llluu] \arrow[rrruuu] \arrow[rrru] &  &                  &                  
\end{tikzcd}
\\

$Y_7:$
\begin{tikzcd}[row sep=.5em]
             & 1 \arrow[rdd] \arrow[rddd] &                           &                            &                           & 5 \arrow[ldd] \arrow[lddd] &               \\
2 \arrow[ru] & 8 \arrow[rdd] \arrow[rd]   &                           & 4 \arrow[llu] \arrow[rru]  &                           & 12 \arrow[ldd] \arrow[ld]  & 7 \arrow[lu]  \\
9 \arrow[ru] &                            & 3 \arrow[llu] \arrow[ru]  & 11 \arrow[llu] \arrow[rru] & 6 \arrow[lu] \arrow[rru]  &                            & 14 \arrow[lu] \\
             &                            & 10 \arrow[llu] \arrow[ru] &                            & 13 \arrow[lu] \arrow[rru] &                            &              
\end{tikzcd}
    \caption{The folded quivers $Y_5$ and $Y_7$}
    \label{fig:x_5}
\end{figure}
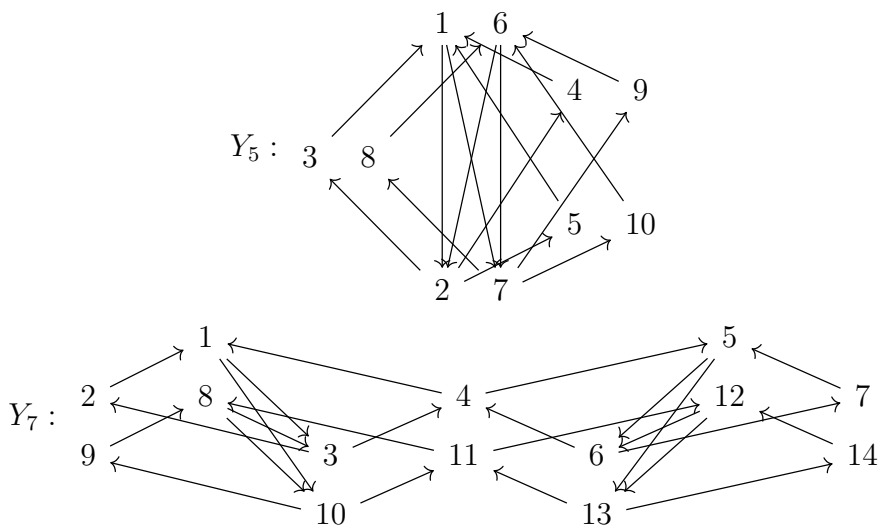

We may also construct two more finite mutation type special foldings from the quiver $X_7$. First by folding all three of the double arrows together, dualizeing and unfolding to produce $Y_5$, shown in figure \ref{fig:x_5}. Then we can repeat this procedure with the $4$ and $5$ groups of $Y_5$ to produce $Y_7$.

\begin{table}[]
    \centering
    \begin{tabular}{c|c}
        Type & size   \\ \hline
        $\widehat{E_6}^1$ & 42 \\
        $\widehat{E_6}^{(1,1)}$ & 55 \\
        $\widehat{E_6}^{(0,2)}$ & 124 \\
        $\widehat{E_6}^{(1,2)}$ & 110 \\
        $\widehat{E_6}^{(2,2)}$ & 23 \\
        
    \end{tabular}
    \quad
        \begin{tabular}{c|c}
        Type & size   \\ \hline

        $\widehat{E_7^1}$ & 192 \\

        $\widehat{E_7}^{(1,1)}$ & 254 \\
        $\widehat{E_7}^{(1,2)}$ & 51 \\

    \end{tabular}
    \quad
        \begin{tabular}{c|c}
        Type &  size   \\ \hline
        $\widehat{E_8^1}$ & 1260 \\
        $\widehat{E_8}^{(1,1)}$ & 1023 \\
        $Y_5$ & 7 \\
        $Y_7$ & 99
        
    \end{tabular}
    \caption{Mutation class sizes of various special foldings}
    \label{tab:my_label}
\end{table}

\section{Punctured Surfaces}\label{sec:punctured}

Our next goal is to construct a folded cluster algebra associated to a  punctured surface where the group mutations correspond to flips that allow for folded triangles.  This will give a geometric interpretation of the quivers $\widehat{D}_n$ and $\widehat{D}_n^1$. Each of the folded cluster algebras constructed will be finite mutation type, as they come from folding a finite mutation type quiver.

We will use the basic notions of cluster algebras associated with triangulations of surfaces and their connections with the Teichmüller space as described in \cite{Triangulated_1,Triangulated_2}. Let $S$ be a marked surface and $\Delta$ an ideal triangulation of $S$
Recall the usual cluster algebra associated with a punctured surface. We denote this cluster algebra $\A(S)$.

Each cluster variable in $\A(S)$ corresponds to an arc on the the surface between two marked points or punctures which may be tagged or untagged according the the rules in \cite{Triangulated_1}. There are however no cluster variables corresponding to arcs that travel from one marked point around a puncture and back to the same marked point.  We denote such an arc by $\l_m^p$ and we call the triangle they cut out a self folded triangle.

We consider a new set of rules for flips of a triangulation on a punctured surface, which we call special flip rules:
\begin{enumerate}
    \item Flips in a punctured digon give a self folded triangle.
    \item The flip of an arc in the interior of a self folded triangle will become a tagged arc.
    \item All tags at a puncture must agree.
\end{enumerate}
These rules completely determine the results of all possible special flips. Moreover, every triangulation of a punctured surface achieved from a sequence of special flips is simply a usual triangulation where each puncture is either fully tagged or untagged. If it possible to place an arc on the surface which cuts out a given puncture as a self folded triangle, then it is possible to achieve triangulations where that puncture has been tagged.

\begin{figure}
    \centering
    \includegraphics[scale=.20]{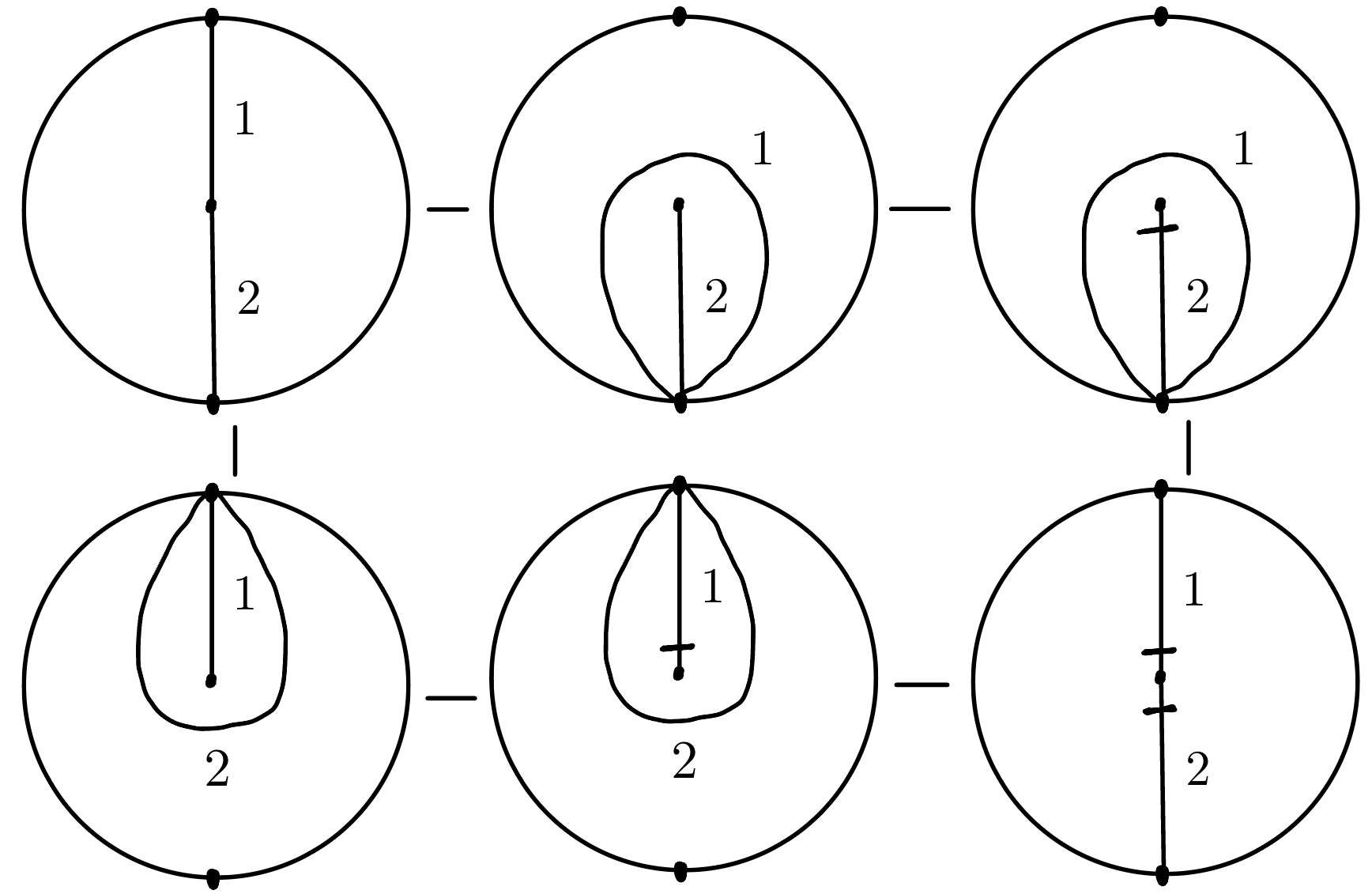}
    \caption{Special flip rules for a punctured digon}
    \label{fig:special_flip_rules}
\end{figure}

These rules are demonstrated for a punctured digon in figure \ref{fig:special_flip_rules}

\begin{remark}
There is no reason we cant pick some punctures to use special flip rules at and some punctures to use the usual tagging rules. Our construction will construct special cluster algebras for these surfaces too. 
\end{remark}

\begin{definition}
We will call a cluster algebra where the mutation combinatorics corresponds to special flip rules a  \emph{special cluster algebra associated to a surface}. We denote this cluster algebra $\As(S)$
\end{definition}
 
We wish to construct such a cluster algebra for every surface with some special punctures punctures Unfortunately, we cannot seem to do so if $S$ is a punctured sphere with only one special puncture, so we exclude this case from the following theorem.

\begin{theorem}
There is a special cluster algebra for any punctured surface with $\chi(S)<0$. Moreover, this cluster algebra is a folded cluster algebra of a cluster algebra associated to a triangulated surface with no special punctures.  
\end{theorem}

Given $S$ a punctured surface with some special punctures, we can compare the usual cluster algebra and the special cluster algebra associated with $S$. Consider a triangulation $\Delta$ of $S$ which contains no tagged arcs or self-folded arcs, and consider the algebras $\A(S)$ and $\As(S)$ both starting with the seed corresponding to $\Delta$

\begin{theorem}
The cluster variable for a tagged arc at a special puncture in $\As(S)$ is exactly 2 times the corresponding cluster variable in $\A(S)$. The cluster variables of $\As(S)$ corresponding to loops $\l_m^p$ are equal to the product of the cluster variables of $\A(S)$ correspond to the tagged and untagged at $p$ arcs from $m$ to $p$. 
\end{theorem}

                                                                                    
Our general strategy for a proof will be to construct a surface, $\Tilde{S}$ which is a two (or four fold cover, in some cases) of $S$ ramified along the special punctures. We consider $\tilde{S}$ as a punctured surface, but with non-special punctures at the ramification points. Then we consider the deck transformation(s) of this cover, $\sigma, (\tau)$, fixing all of ramification points and a regular triangulation where every puncture is either untagged or fully tagged of $\Tilde{S}$ which is invariant under $\sigma$. Such a triangulation gives a triangulation of $S$ with possibly self folded triangles.
The folding of the quiver associated to the triangulation of $\tilde{S}$ will give the folded seed for our special cluster algebra associated to $S$.

This folded seed is valid and satisfies the folded cluster condition by our lemmas \ref{thm:order_2folds},\ref{thm:folded_clusters}; in the case of a four fold cover, this follows since this cover factors as two two-fold covers. 

The group mutations of this folding clearly correspond to our special flip rules; The flips either are corresponding to the usual flip rules or are happening in a special punctured digon on $S$. Each special punctured digon on $S$ lifts to a regular punctured square on $\tilde{S}$ where one easily checks that the group mutations give the special flip rules. 

The facts about the values of the cluster variables by direct calculation of the cluster variables on $\tilde{S}$. Since we start with a collection of cluster variables which are symmetric under the deck transformations, the associated lambda lengths parameterize a point in Teichmüller space with a symmetric hyperbolic metric. Since the folding is valid and cluster we know that every group mutation produces a new symmetric triangulation, and so the new lambda lengths will also be symmetric. 

We essentially follow the construction in section 12 of \cite{Felikson_triangulated_orbifolds} to construct $\Tilde{S}$, but we reiterate it here in the next subsections for completeness. The issue with a sphere with one special puncture also appears in their construction.


\subsection{Punctured disk}

First we construct a covering special covering for a multi punctured monogon, with all special punctures. This construction depends on the parity of the number of special punctures, $p$. If $p=2g+1$ is odd, then we take $\tilde{S}$ to be a surface of genus $g$ with $p$ punctures and one boundary component with two marked points. The covering is given by the hyperelliptic involution fixing the $p$ punctures and and the boundary component (swapping the two marked points). 

If $p=2g+2$, then we take $\tilde{S}$ to be the surface of genus $g$ with $p$ punctures and two boundary components with one marked point each, and the covering again given by the hyperelliptic involution fixing all the punctures, but swapping the two boundary components. 

Each of these surfaces has a usual triangulation which is invariant under the covering involution. Note that in either case, $\tilde{S}$ has two boundary arcs, either on two separate components or on the same component. 

The cases where $S$ has one or two special punctures give the $\Dhat{n}$ and $\widehat{D^1_{n+3}}$ quivers respectively.

\subsection{Surfaces with a boundary component}

We now assume that $S$ is a punctured surface with $n$ punctures surface with at least one boundary component with at least one marked point. Choose one such marked point and call it $m$. Let $l$ be an arc from $m$ surrounding all of the special punctures and returning to $m$. Thus, inside $l$ we have a $p$ special punctured monogon and outside $l$ there is a surface with no special punctures. Call this punctured monogon $M$.

We will build the surface $\tilde{S}$ along with its involution for this surface by gluing the surface $\tilde{M}$ previously constructed with two copies of $S-M$ along the two boundary arcs of $\tilde{M}$. The hyperelliptic involution on $\tilde{M}$ extends to an involution on $\tilde{S}$ by swapping the two copies. 

We can use any symmetric triangulation of $\tilde{M}$ along with any triangulation of $S-M$ to make a symmetric triangulation of $\tilde{S}$. 

\begin{figure}
    \centering
    \begin{subfigure}{.45\textwidth}
    \includegraphics[scale=.1]{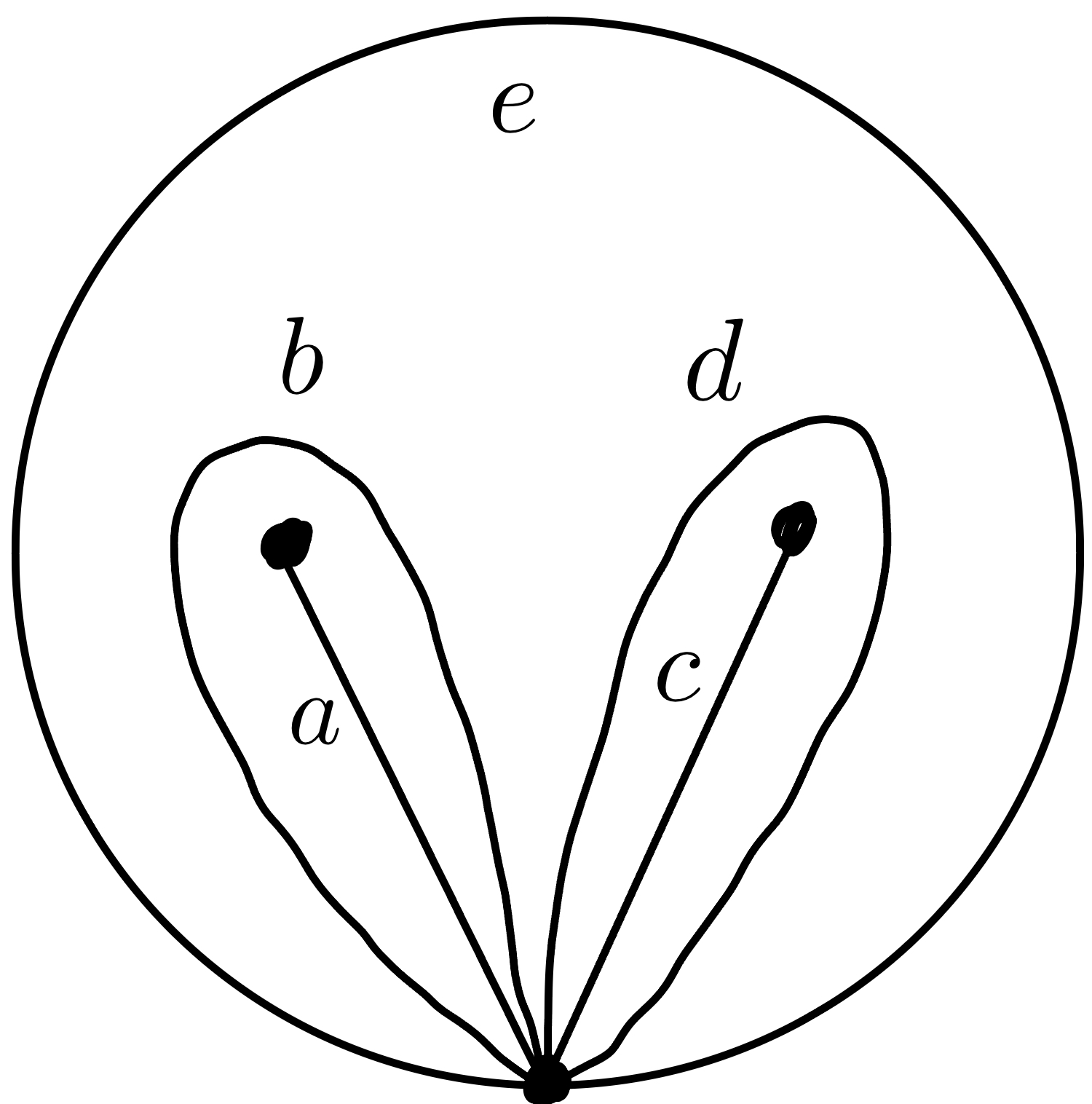}
    \caption{A triangulation of a twice punctured digon.}
    \end{subfigure}
    \begin{subfigure}{.45\textwidth}
    \includegraphics[scale=.15]{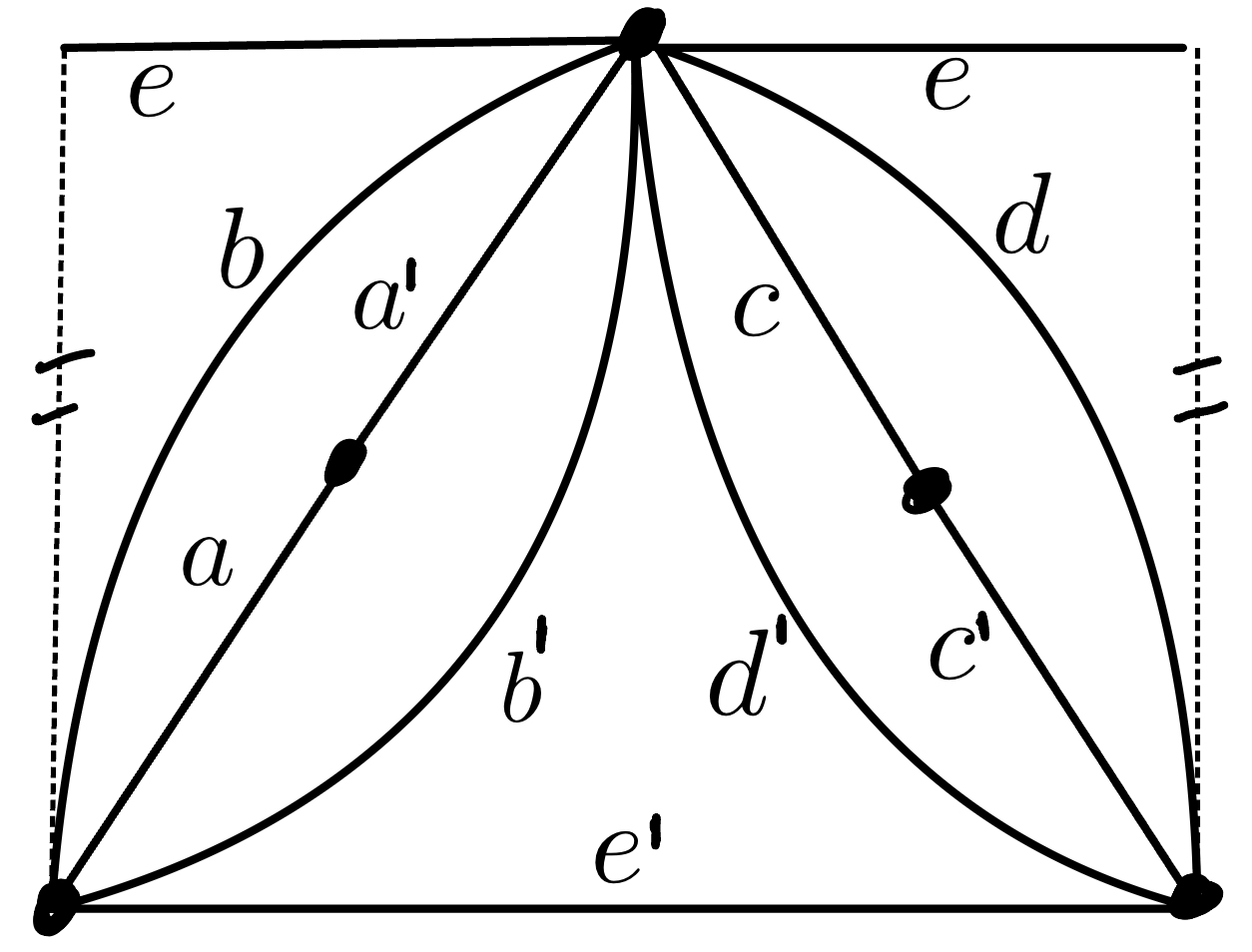}
    \caption{The two fold cover.}
    \end{subfigure}
    \label{fig:twice_puncutred}
\end{figure}





\subsection{Closed surfaces with punctures}



Now consider $S$ a closed surface with $p$ special punctures. First we assume $p>1$, even. In this case we construct $\tilde{S}$ directly by pairing up the special punctures into groups $p_i,p_{i+1}$. Then we cut $S$ along simple arcs connecting each pair creating one boundary with two marked points for each pair. We glue two copies of the resulting surface along each boundary component to make $\tilde{S}$. The involution is simply the involution fixing all of the special punctures and swapping the two copies. We can easily make a symmetric triangulation of $\tilde{S}$ by starting with any triangulation of $S$ which contains all of the arcs we cut along. 

If $p>1$ is odd, we first pair up all of the punctures except one, call it $p_0$. Cutting $S$ along each pair and doubling results in a surface $S_1$ with an involution fixing all of the special punctures except $p_0$, which is mapped to $p_0'$. Now we consider $S_1$ as a closed surface with two special punctures, $p_0$ and $p_0'$. Repeating this procedure will give $\tilde{S}$ as a four fold cover of $S$ ramified over all of the special punctures.

Finally if $p=1$ we use the same trick as \cite{Felikson_triangulated_orbifolds}. We assume that $S$ is not a sphere so that we may find an arc from $p$ to itself that does not separate the surface into two pieces. We cut $S$ along this arc, and glue two copies of the resulting surface along their boundaries to create $S_1$ a surface which is a two fold unramified cover of $S$. $S_1$ now has 2 special punctures so be can construct $\tilde{S}$ from it with a four fold cover over $S$ ramified at $p$.

\subsection{The algebra $\hat{D}_{n}$}\label{sef:dhat_disk}

The folding $\Dhat{n}$ is associated to a disk with $n$ marked points and one special puncture. We can use this correspondence to count the total number of cluster variables and the total number of clusters. 

The number of cluster variables is given by $n(n+1)$ since the number of arcs on a punctured disk with special flip rules is equal to the number of arcs coming from the normal flip rules plus the number of arcs forming self folded triangles. Thus there are $n^2$ arcs coming from the cluster variables in the usual $D_n$ cluster algebra plus $n$ arcs corresponding to the folded triangles. 

The number of clusters is also easy to compute. They can be organized into disjoint sets as follows:
\begin{enumerate}
    \item Triangulations with no folded triangles
    \item Triangulations with a folded triangle containing an untagged arc
    \item Triangulations with a folded triangle containing a tagged arc
\end{enumerate}
The total of the first two sets is equal to the number of triangulations of a punctured disk with a normal puncture, which can be seen by replacing the folded arc with a tagged arc to the puncture. The size of the third set is equal to the number of triangulations of an $n+1$-gon times $n$. 
Thus the total is given by $\frac{3n-2}{n}\binom{2n-2}{n-1} + \frac{n}{n}\binom{2n-2}{n-1} = \binom{2n}{n} $

\begin{remark}
As remarked earlier, this result echos the calculation for the $B_n$ and $C_n$ cluster algebras which are associated to an $n+1$-gon with one orbifold point. However the surfaces and their symmetries are different: a punctured $n$-gon vs an $n+1$-gon with an orbifold point.
\end{remark}

\section{Permutohedra}\label{sec:permuta}

We can construct a special folded cluster algebras for which the exchange complex is a generalized permutohedra for each Dynkin type. The construction is simple enough. 

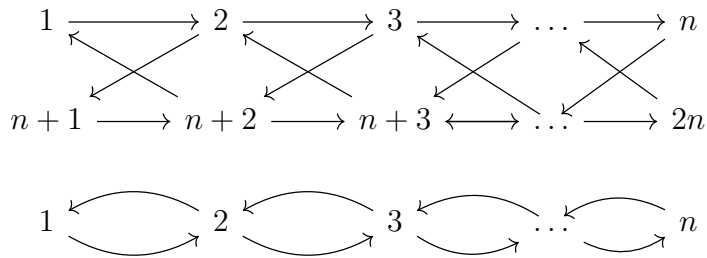
\begin{figure}
    \centering
\begin{tikzcd}
1 \arrow[r]             & 2 \arrow[ld] \arrow[r]                        & 3 \arrow[ld] \arrow[r]                                 & \dots  \arrow[r] \arrow[ld]              & n \arrow[ld]            \\
n+1 \arrow[r]           & n+2 \arrow[lu] \arrow[r]                      & n+3 \arrow[lu] \arrow[r]                                & \dots \arrow[l] \arrow[r] \arrow[lu]              & 2n \arrow[lu]           \\
1 \arrow[r, bend right] & 2 \arrow[l, bend right] \arrow[r, bend right] & 3 \arrow[l, bend right] \arrow[r, bend right] & \dots \arrow[l, bend right] \arrow[r, bend right] & n \arrow[l, bend right]
\end{tikzcd}
    \caption{The quiver $\Bar{A_n}$}
    \label{fig:A_n_permutohedron}
\end{figure}

Let $X_n$ be an ADE type Dynkin diagram. We define a quiver $\Bar{X_n}$ with $2n$ and a folding with $n$ groups by replacing each edge of $X_n$ with a two-cycle edge, see figure \ref{fig:A_n_permutohedron} for the $\Bar{A_n}$ case. 

\begin{theorem}
$\Bar{X_n}$ satisfies the folding condition and has an associated special folded cluster algebra. The exchange complex of this special cluster algebra is the permutohedron of type $X_n$.
\end{theorem}
\begin{proof}

One easily checks that every group mutation followed by swapping the two members of the group preserves the quiver isomorphism class. Thus these quivers satisfy the folding condition. Moreover, it is also clear that there is a folded cluster algebra associated to each quiver. It just remains to be seen that the exchange complex is the corresponding generalized permutohedron. 

We can see that the cluster modular group of $\Bar{X_n}$ is at least a quotient of $W(X_n) \rtimes Aut(X_n)$ in each case in the following way: Since each mutation preserves the quiver isomorphism class, the cluster modular group is generated by the mutations $\mu_i$ for each group along with the elements of $Aut(X_n)$. The mutations clearly satisfy the Coxeter relations since the only relations needed to be checked are those coming from $A_2$, $B_2$ and $G_2$ which can be checked by hand.

We only need to see now that the cluster modular group is not a smaller quotient of the corresponding Coxeter group. Following \cite{Maxwell_normal_subgroups}, there is only a finite list of possible quotients of Coxeter groups and we can see easily that none of these are possible since each quotient listed identifies some generators of the Coxeter group, but each generator of the cluster modular group is clearly different. 

\end{proof}

Some of these quivers are in known mutation classes.
The unfolding of the $\Bar{A_3}$ quiver is mutation equivalent to that of a four punctured sphere (with regular punctures). This folding gives the special cluster algebra for a three punctured sphere with two special punctures.

The unfolding of $\Bar{A_5}$ is mutation equivalent to the unfolding of $Y_5$.  The unfolding of $\Bar{D_4}$ is the cluster type of the Fock-Goncharov moduli space associated to $SL_3$ on a three punctured sphere.

\bibliography{References}

\appendix
\section{Exchange Graph Pictures}\label{sec:appendix}

In the following pictures, the nodes represent clusters, the edges mutations, and the colors approximate the quiver isomorphism class of the quiver underlying each cluster.

\begin{figure}[b]
    \centering
    \begin{subfigure}{.48\textwidth}
    \centering
    \includegraphics[scale=.17]{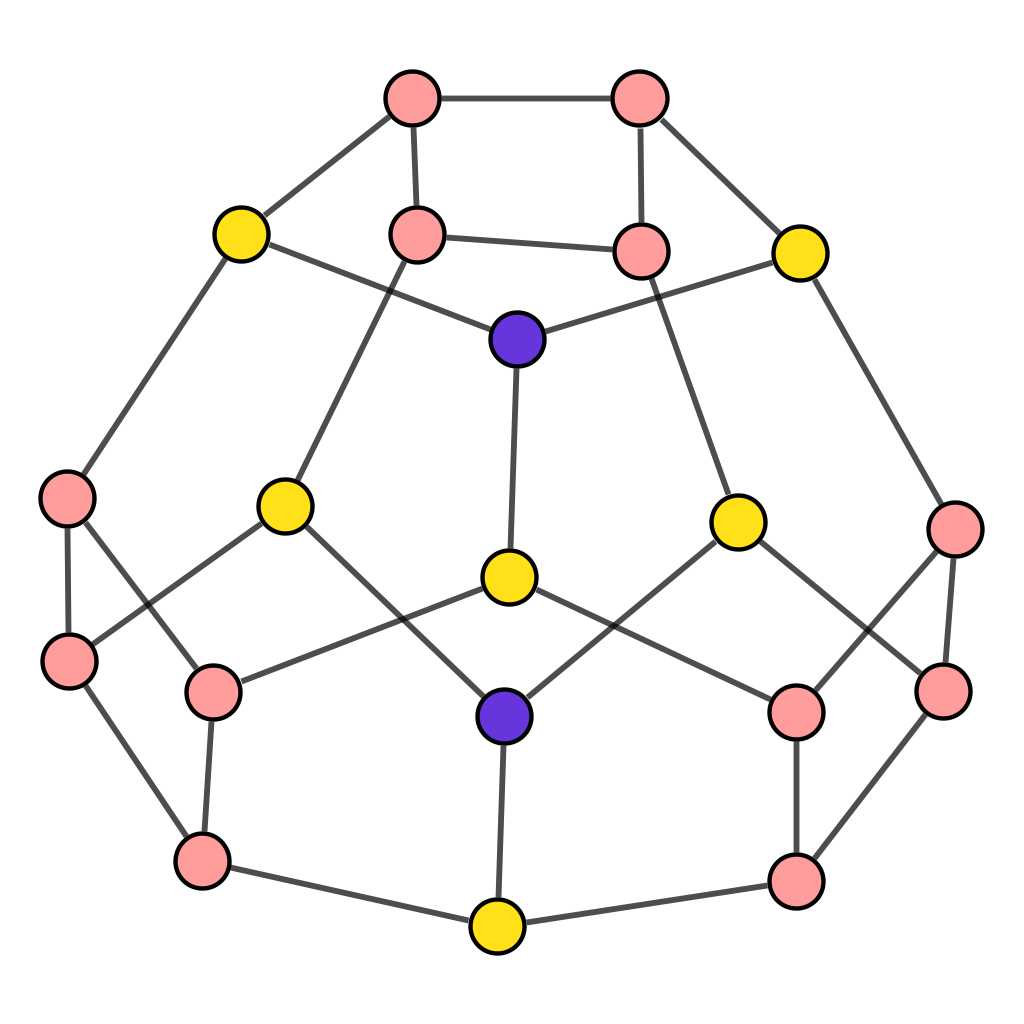}
    \caption{$\widehat{D_3}$}
    \end{subfigure}
    \begin{subfigure}{.48\textwidth}
    \centering
    \includegraphics[scale=.17]{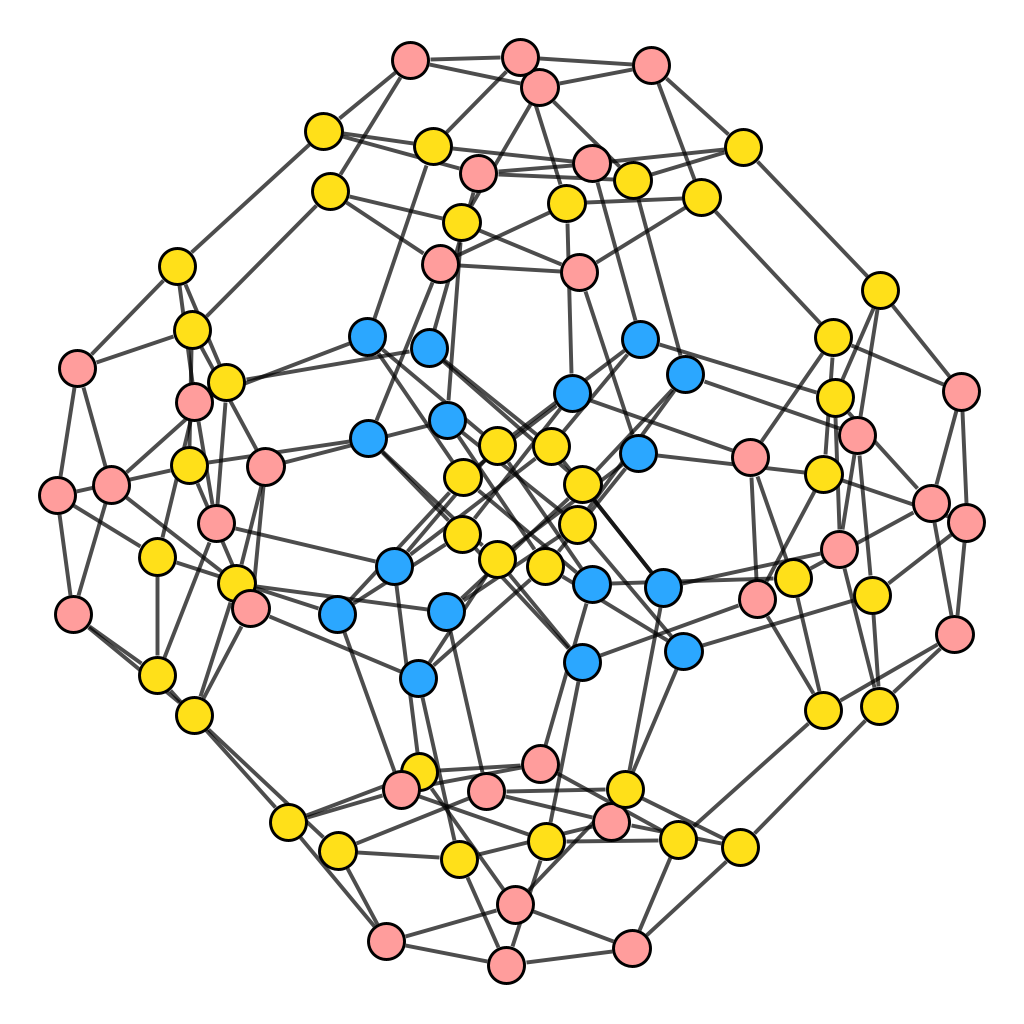}
    \caption{$\widehat{E_6}$}
    \end{subfigure}
\end{figure}
\begin{figure}
\centering
    \begin{subfigure}{.48\textwidth}
    \centering
    \includegraphics[scale=.16]{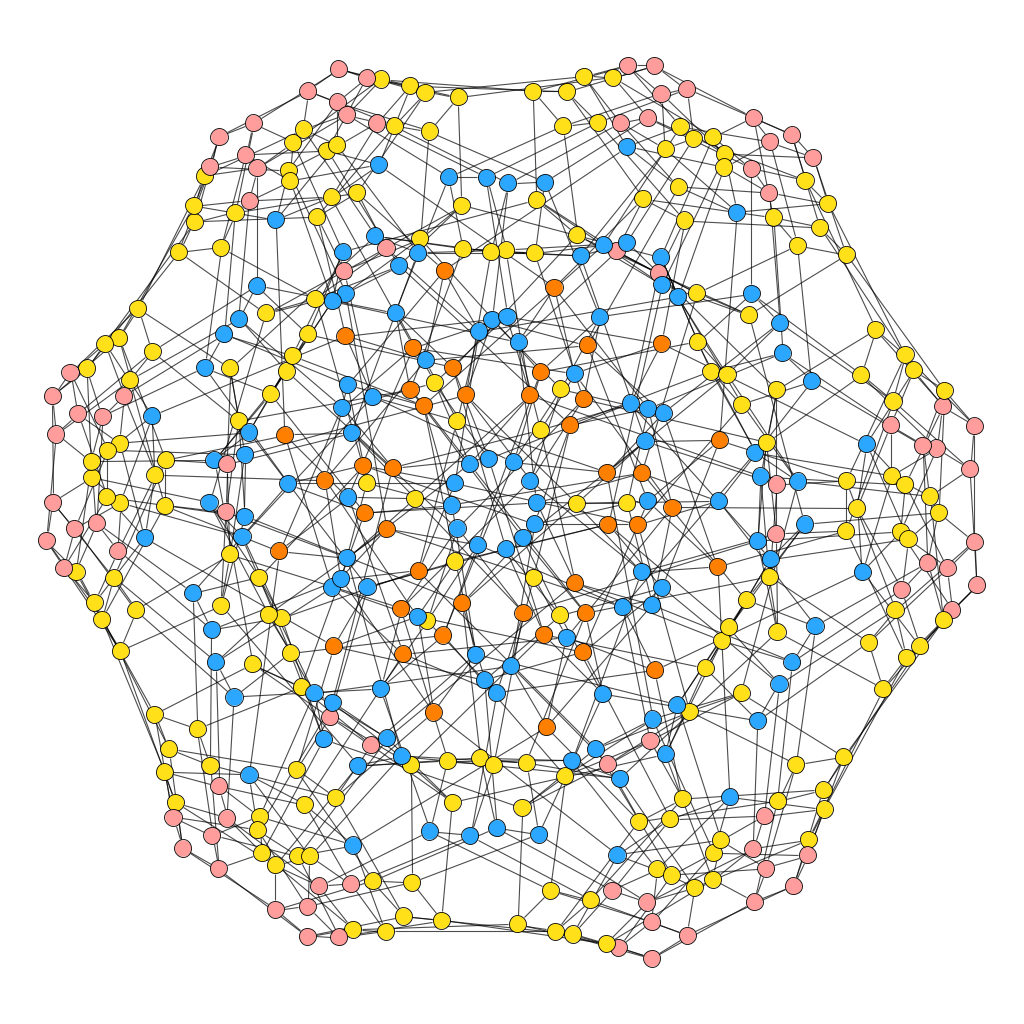}
    \caption{$\widehat{E_7}$}
    \end{subfigure}
    \begin{subfigure}{.48\textwidth}
    \centering
    \includegraphics[scale=.16]{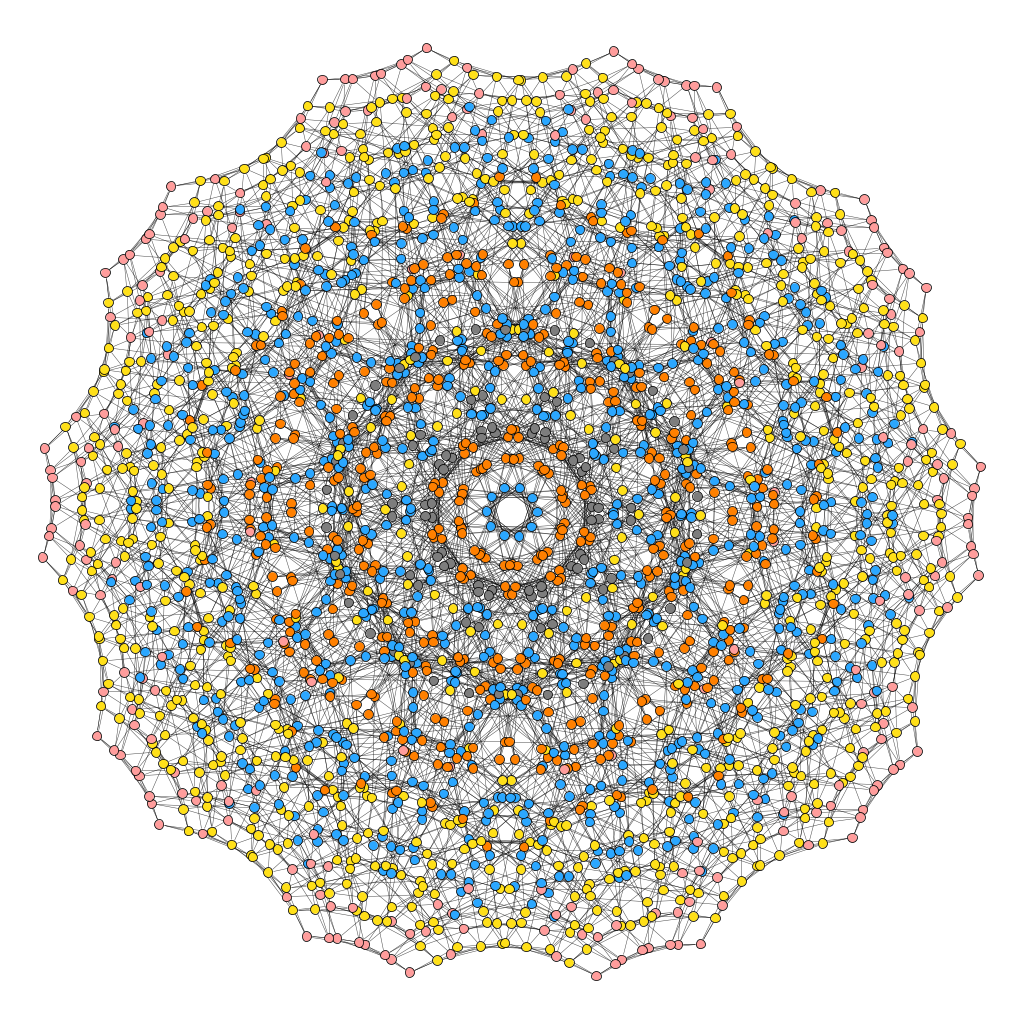}
    \caption{$\widehat{E_8}$}
    \end{subfigure}
    \begin{subfigure}{.48\textwidth}
    \centering
    \includegraphics[scale=.15]{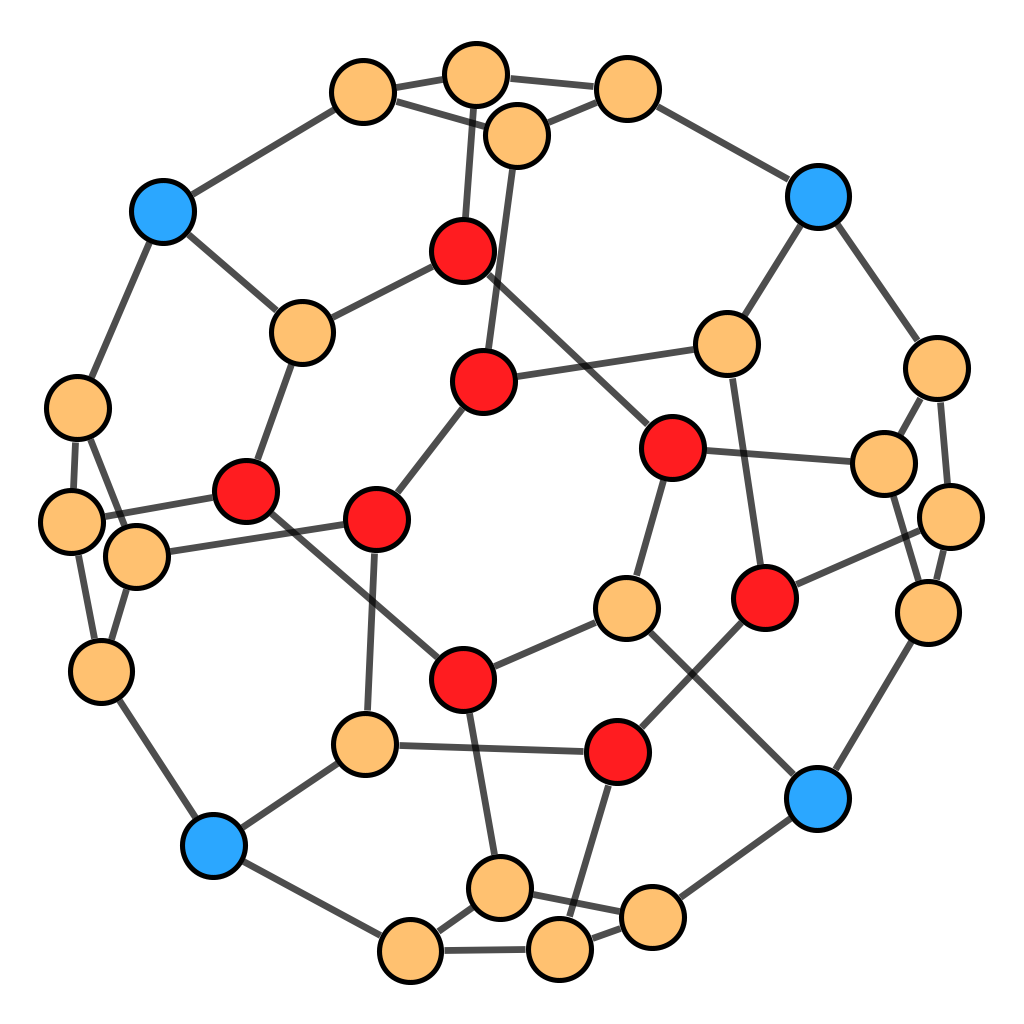}
    \caption{$\widehat{H_3}$}
    \end{subfigure}
    \begin{subfigure}{.48\textwidth}
    \centering
    \includegraphics[scale=.15]{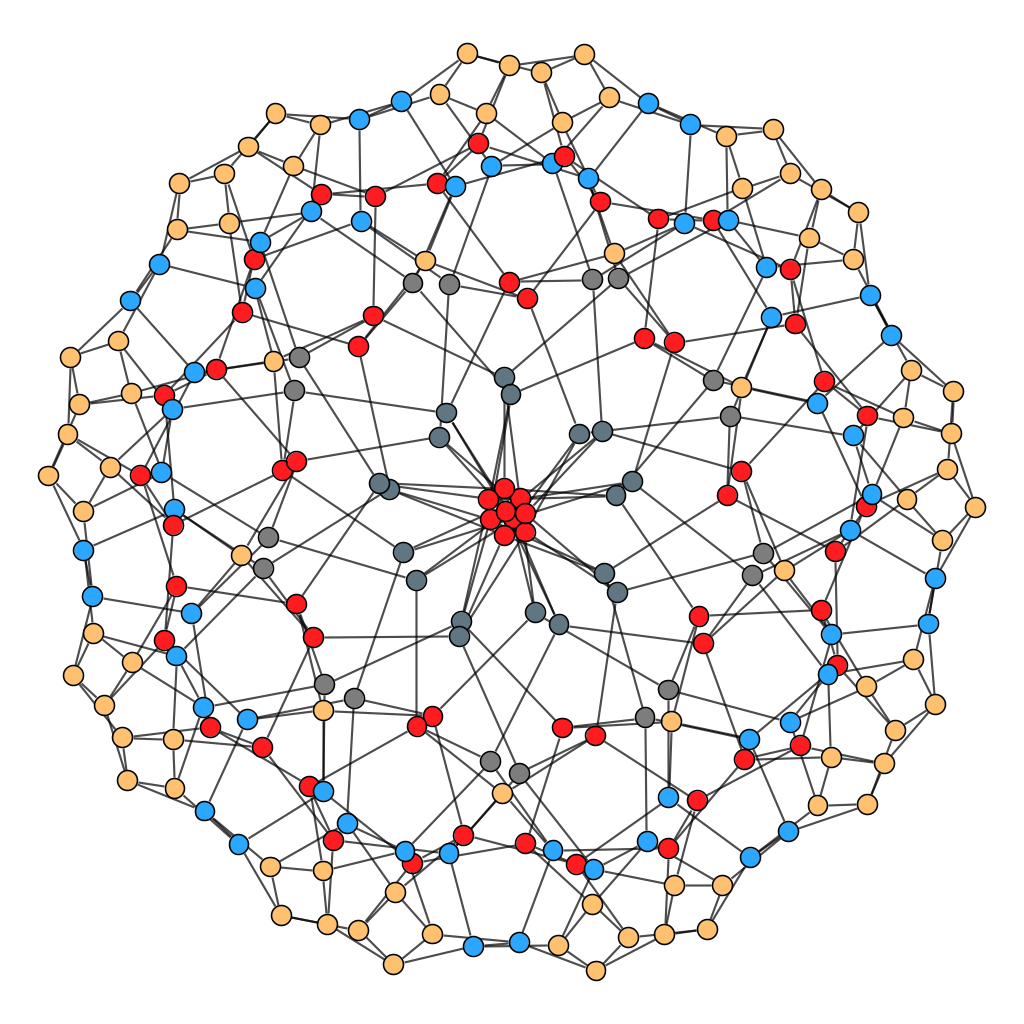}
    \caption{$\widehat{H_4}$}
    \end{subfigure}
    \label{fig:exchange graphs}
\end{figure}
\begin{figure}
        \centering
    \includegraphics[scale=.29]{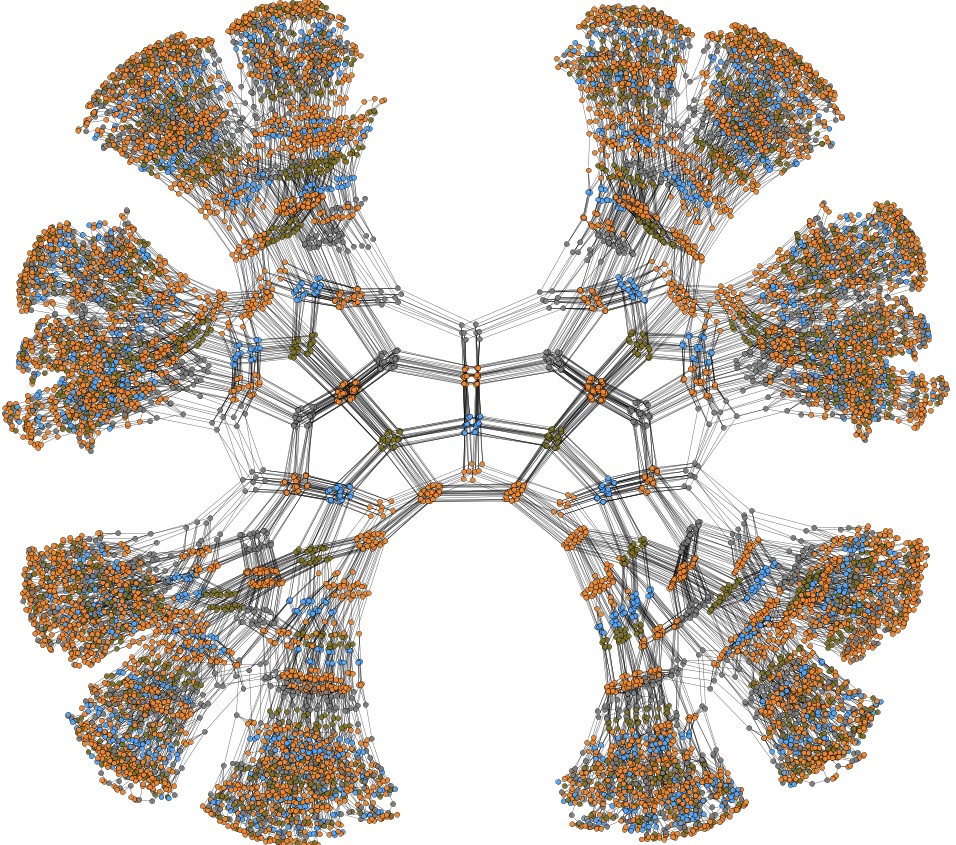}
    \caption{$\widehat{Y_5}$}
    \end{figure}
\begin{figure}
\centering
    
    \begin{subfigure}{.95\textwidth}
    \centering
    \includegraphics[scale=.4]{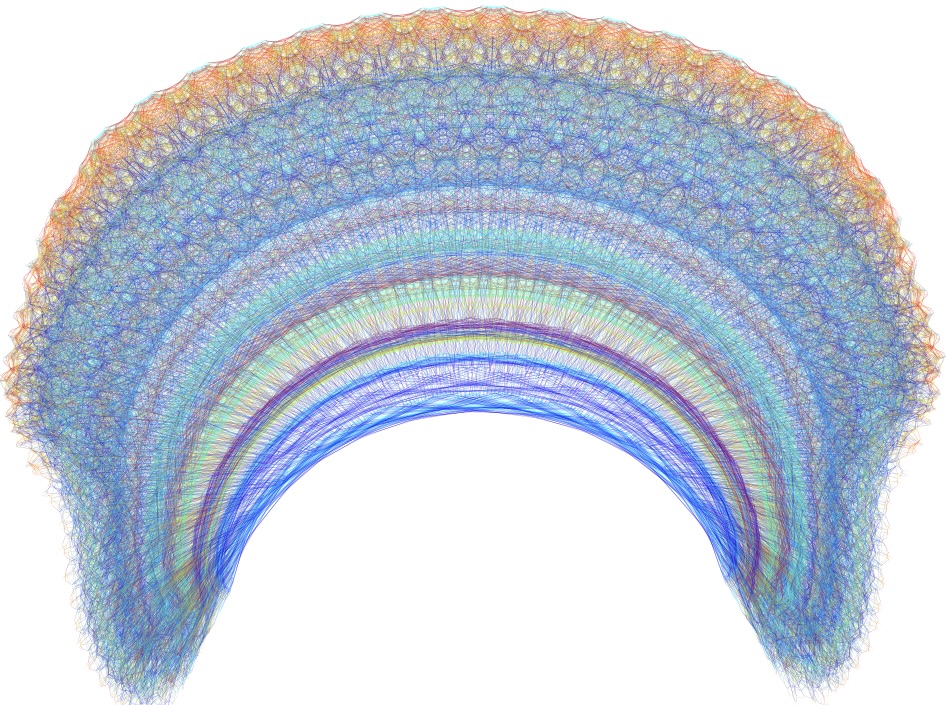}
    \caption{$\widehat{E_8}^1$}
    \end{subfigure}
    \begin{subfigure}{.95\textwidth}
    \centering
    \includegraphics[scale=.4]{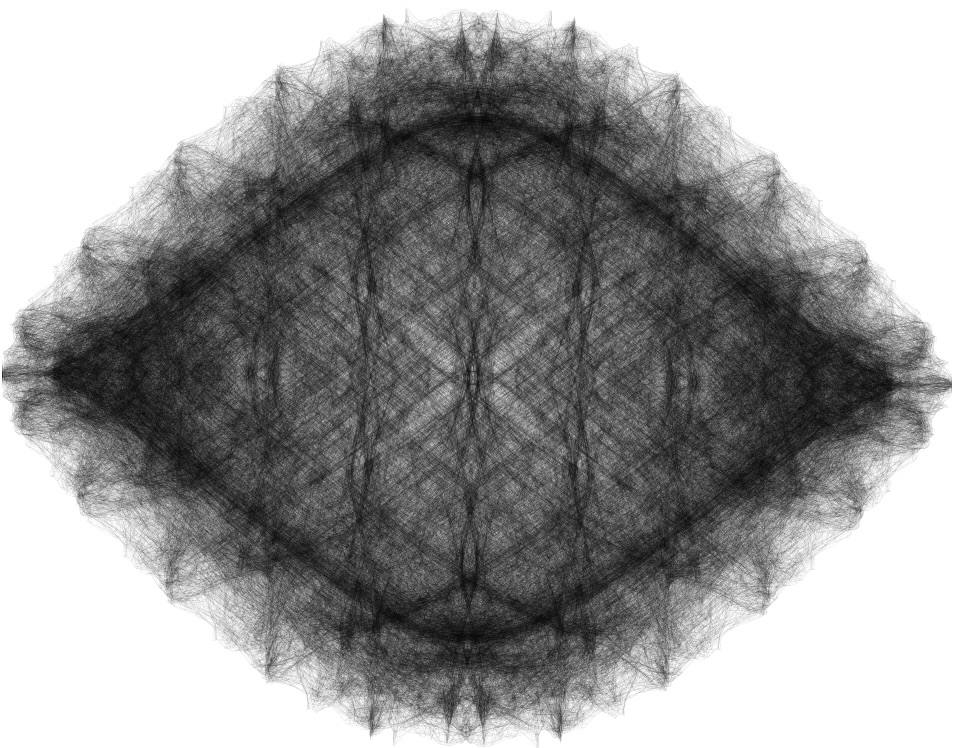}
    \caption{$\widehat{E_7}^{(1,2)}$}
    \end{subfigure}
    \label{fig:my_label}
\end{figure}

\end{document}